\documentclass{article}
\usepackage{latexsym,amsmath,amssymb,amscd}

\begin{document}

\title{On spherically symmetric solutions of the Einstein-Euler equations}
\author{Tetu Makino \footnote{Department of Applied Mathematics, Yamaguchi University, Ube, 755-8611 Japan, E-mail: makino@yamaguchi-u.ac.jp}}
\date{\today}
\maketitle

\newtheorem{Lemma}{Lemma}
\newtheorem{Proposition}{Proposition}
\newtheorem{Theorem}{Theorem}
\newtheorem{Definition}{Definition}

\numberwithin{equation}{section}

\begin{abstract}
We construct spherically symmetric solutions to the Einstein-Euler equations, which give models of gaseous stars in the framework of the  general theory of relativity. We assume a realistic barotropic equation of state. Equilibria of the spherically symmetric Einstein-Euler equations are given by the Tolman-Oppenheimer-Volkoff equations, and time periodic solutions around the equilibrium of the linearized equations can be considered. Our aim is to find true solutions near these time-periodic approximations. Solutions satisfying so called physical boundary condition at the free boundary with the vacuum will be constructed using the Nash-Moser theorem. This work also can be considered as a touchstone in order to estimate the universality of the method which was originally developed for the non-relativistic Euler-Poisson equations.
\end{abstract}

{\it Key Words and Phrases.} Einstein equations, Spherically symmetric solutions, Vacuum boundary, Nash-Moser theorem

{\it 2010 Mathematical Subject Classification Numbers.} 35L05, 35L52, 35L57, 35L70, 76L10, 76N15, 83C05, 85A30

\section{Introduction}
Recently U. Brauer and L. Karp \cite[Theorem 2.3]{BrauerK} established a local existence theorem of solutions to the Cauchy problem for the Einstein-Euler equations, which describes a relativistic self-gravitating perfect fluid having density either compactly supported or falling off at infinity in an appropriate manner.

 In their work \cite{BrauerK} the energy-momentum tensor of the perfect fluid takes the form
$$T^{\mu\nu}=(\epsilon +P)U^{\mu}U^{\nu}-Pg^{\mu\nu},$$
where $\epsilon=c^2\rho$ is the energy density, $P$ is the pressure, and $U^{\mu}$ is the velocity 4-vector. Here it is assumed that $P=K\epsilon^{\gamma}, K>0, \gamma>1$, and the quantity $$w:=\epsilon^{\frac{\gamma-1}{2}}=c^{\gamma-1}
\rho^{\frac{\gamma-1}{2}}$$
is introduced. The main result requires that the initial data satisfy $w\in H_{s+1}$ with $s>3/2$ so that $w\in C^1$ at least. 

However a spherically symmetric equilibrium, which solves the
Tolman-\\
Oppenheimer-Volkoff equation, satisfies $w \sim \mbox{const}
(r_+-r)^{1/2}$ as
$r\rightarrow r_+-0$ provided that the equilibrium has a finite radius $r_+$.  
See Section 3. Hence such an equilibrium is excluded from the class of density distributions admissible to this local existence theorem. We are faced with the same situation in the non-relativistic problem governed by the Euler-Poisson
equations as discussed in \cite{Makino1986}.

Recently this trouble was partially overcome by \cite{MakinoOJM} in the Euler-Poisson equations for the non-relativistic case. So, a similar discussion is required for the relativistic problem. It is the aim of this article.

\section{Spherically symmetric evolution equations}
The Einstein equations read (\cite[\!(95.5)\!]{LandauL}.)
\begin{equation}
R_{\mu\nu}-\frac{1}{2}g_{\mu\nu}R=\frac{8\pi G}{c^4}T_{\mu\nu}.
\end{equation}
Here $R_{\mu\nu}$ is the Ricci tensor and $R$ is the scalar curvature $g^{\alpha\beta}R_{\alpha\beta}$ associated with the metric
\begin{equation}
ds^2=g_{\mu\nu}dx^{\mu}dx^{\nu}, \label{ds2}
\end{equation}
and $T^{\mu\nu}$ is the energy-momentum tensor of the matter. $G$ is the constant of gravitation (6.67$\times 10^{-8} \mbox{cm}^3/\mbox{g}\cdot \mbox{sec}^2$), and $c$ is the speed of light (3.00$\times 10^{10}\mbox{cm}/\mbox{sec}$). The Einstein equations (1) imply the Euler equations
\begin{equation}
\nabla_{\nu}T^{\mu\nu}=0,
\end{equation}
where $\nabla$ denotes the covariant derivative associated with the metric (\ref{ds2}). The details can be found in \cite{LandauL} or \cite{MisnerTW}. 

The energy-momentum tensor of a perfect fluid is given by 
\begin{equation}
T^{\mu\nu}=(c^2\rho +P)U^{\mu}U^{\nu}-Pg^{\mu\nu},
\end{equation}
(\cite[94.4)\!]{LandauL}),
where $\rho$ is the mass density, $P$ is the pressure, and $U^{\mu}$ stands for the 4-dimensional velocity vector such that
$U^{\mu}U_{\mu}=1$. In this article we always assume\\

(A0)\hspace{10mm} {\bf $P$ is a given analytic function of $\rho>0$ such that $0<P, 0<dP/d\rho <c^2$ for $\rho >0$ and $P\rightarrow 0$ as $\rho \rightarrow +0$.} \\

If we assume the spherical symmetry, the Einstein-Euler equations are reduced as follows.

We consider the metric of the form
\begin{equation}
ds^2=e^{2F}c^2dt^2-e^{2H}dr^2-
R^2(d\theta^2+\sin^2\theta d\phi^2), \label{metric}
\end{equation}
(\cite[p.304, (1)\!]{LandauL}),
where $F, H$ and $R$ are functions of $t, r(\geq 0)$. (Here $R$ does not mean the scalar curvature $g^{\mu\nu}R_{\mu\nu}$.) Then the non-zero components of the Einstein tensor $G_{\mu}^{\nu}:=R_{\mu}^{\nu}-
\frac{1}{2}\delta_{\mu}^{\nu}R$, where $R$ is the scalar curvature,  are (see \cite[p.305, (2)(3)(4)(5)\!]{LandauL}):
\begin{align*}
G_0^0&=e^{-2H}\Big(-\frac{R'^2}{R^2}-2\frac{R''}{R}+
2\frac{H'R'}{R}\Big)+
e^{-2F}\Big(\frac{\dot{R}^2}{R^2}+2\frac{\dot{H}\dot{R}}{R}\Big)+\frac{1}{R^2}, \\
G_1^1&=e^{-2F}\Big(\frac{\dot{R}^2}{R^2}+
2\frac{\ddot{R}}{R}-2\frac{\dot{F}\dot{R}}{R}\Big)
-e^{-2H}\Big(\frac{R'^2}{R^2}+2\frac{F'R'}{R}\Big)+\frac{1}{R^2}, \\
G_2^2&=G_3^3=
e^{-2H}\Big(-\frac{R''}{R}-F''-F'^2+H'F'+\frac{H'R'}{R}-\frac{F'R'}{R}\Big)+ \\
&+e^{-2F}\Big(\frac{\ddot{R}}{R}+\ddot{H}+\dot{H}^2-\dot{H}\dot{F}
+\frac{\dot{H}\dot{R}}{R}-\frac{\dot{F}\dot{R}}{R}\Big), \\
e^{2H}G_0^1&=-e^{2F}G_1^0=
2\Big(\frac{\dot{R}'}{R}-\frac{\dot{H}R'}{R}-\frac{F'\dot{R}}{R}\Big).
\end{align*}
Here $\dot{A}$ stands for $\partial A/c\partial t$ and $A'$ stands for
$\partial A/\partial r$. 
Of course the coordinates $x^{\mu}$ are taken as
$$x^0=ct,\quad x^1=r,\quad x^2=\theta,\quad x^3=\phi. $$

By a freedom of choice of $r$ we take it in such a way that the flow is apparently static, say, we suppose
\begin{equation}
U^0=e^{-F}, \quad U^1=U^2=U^3=0.
\end{equation}
Then the energy-momentum tensor turns out to be
\begin{equation}
T_0^0=c^2\rho, \quad T_1^1=T_2^2=T_3^3=-P, \quad T_0^1=T_1^0=0.
\end{equation}

The equation $\nabla_{\mu}T_0^{\mu}=0$ gives
\begin{equation}
c^2\dot{\rho}+\Big(\dot{H}+\frac{2\dot{R}}{R}\Big)(c^2\rho +P)=0, \label{23}
\end{equation}
and the equation $\nabla_{\mu}T_1^{\mu}=0$ gives
\begin{equation}
P'+F'(c^2\rho+P)=0. \label{24}
\end{equation}
By integrating (\ref{24}) we can suppose that $F$ is a function of $\rho$ given by
\begin{equation}
F=F(\rho)=-\int^{\rho}
\frac{1}{c^2\rho + P}\frac{dP}{d\rho}d\rho. \label{25}
\end{equation}
Let us introduce the variable $m$ by
\begin{equation}
m=4\pi\int_0^R\rho R^2dR=
4\pi\int_0^r\rho R^2R'dr. \label{26}
\end{equation}
The variable $V$ is defined by
\begin{equation}
V=ce^{-F}\dot{R}. \label{27}
\end{equation}
Then the equation $G_0^1=0$ turns out to be
\begin{equation}
\dot{H}=\frac{1}{c}e^F\frac{V'}{R'}. \label{28}
\end{equation}
Substituting (\ref{27})(\ref{28}) into (\ref{23}), we have
\begin{equation}
c^2\dot{\rho}=-\frac{1}{c}e^F(c^2\rho+P)\Big(\frac{V'}{R'}+\frac{2V}{R}\Big).
\label{29}
\end{equation}
Eliminating the time derivatives from the equation 
$\displaystyle G_0^0=\frac{8\pi G}{c^2}\rho $, we have
$$\frac{8\pi G}{c^2}\rho R^2R'=\Big(
-RR'^2e^{-2H}+\frac{1}{c^2}RV^2+R\Big)'.
$$
Integrating this, keeping in mind that $R$ should vanish at $r=0$, we get
\begin{equation}
m=\frac{c^2R}{2G}\Big(
\frac{V^2}{c^2}+1-R'^2e^{-2H}\Big), \label{30}
\end{equation}
from which we get
\begin{equation}
e^{2H}=\Big(1+\frac{V^2}{c^2}-\frac{2Gm}{c^2R}\Big)^{-1}R'^2. \label{30bis}
\end{equation}
Differentiating (\ref{27}) with respect to $t$ and using the equation
$\displaystyle G_1^1=-\frac{8\pi G}{c^4}P$ and (\ref{30}), we obtain
\begin{equation}
\frac{\dot{V}}{c}e^{-F}=-\frac{GR}{c^2}
\Big(\frac{m}{R^3}+\frac{4\pi P}{c^2}\Big)
-e^{-2H}\frac{R'P'}{c^2\rho+P}, \label{31}
\end{equation}
or, from (\ref{30bis}),
\begin{equation}
e^{-F}c\dot{V}=
-GR\Big(\frac{m}{R^3}+\frac{4\pi P}{c^2}\Big)
-\Big(1+\frac{V^2}{c^2}-\frac{2Gm}{c^2R}\Big)
\frac{P'}{R'(\rho+P/c^2)}. \label{31bis}
\end{equation}
Differentiating (\ref{30}) with respect to $t$ and using the equation
$G_1^0=0$, we have
\begin{equation}
\dot{m}e^{-F}=-\frac{4\pi R^2}{c^3}P V. \label{32}
\end{equation}

Now the equations (\ref{27})(\ref{29})(\ref{31bis})(\ref{32}) govern the evolution of unknowns $R, H, \rho, V, m$. The system of equations to be studied is:
\begin{subequations}
\begin{align}
e^{-F}c\dot{R}&= V  \label{Eq1}\\
e^{-F}c\dot{\rho}&=-(\rho+P/c^2)\Big(\frac{V'}{R'}+\frac{2V}{R}\Big)  
\label{Eq2} \\
e^{-F}c\dot{V}&=-GR\Big(\frac{m}{R^3}+\frac{4\pi P}{c^2}\Big)
-\Big(1+\frac{V^2}{c^2}-\frac{2Gm}{c^2R}\Big)
\frac{P'}{R'(\rho+P/c^2)} \label{Eq3}\\
e^{-F}c\dot{m}&=-\frac{4\pi}{c^2}R^2PV \label{Eq4}
\end{align}
\end{subequations}
Of course we assume (\ref{25}) and (\ref{26}). 
The above equations were derived by 
\cite{MisnerS}. The equations (\ref{Eq1}), (\ref{Eq2}), (\ref{Eq3}),
(\ref{Eq4}) are none other than 
(1.12-$R$), (8.11), (1.12-$U$), (1.12-$m$) of
\cite{MisnerS} respectively. \\

The system of coordinates $(t,r)$ is a co-moving Lagrangian
system of coordinates moving at each point with the fluid. Therefore
if $\rho >0$ for $0\leq r <r_+$ and $\rho=0$ for $r_+\leq r$ at $t=0$, then it remains so for all small $t>0$ along the time  evolution as long as the $C^1$ solution exists, while the surface radius $r_+$ is constant. (Of course the value of 
$R$ at the surface can change in time.) Especially we have $m=m_+$ is constant at $r=r_+$ for all $t>0$. Hence we can take $(t,m)$ as another system of co-moving Lagrangian coordinates. Then we have the formula
\begin{align}
\Big(\frac{\partial}{\partial t}\Big)_r&
=\Big(\frac{\partial}{\partial t}\Big)_m-\frac{4\pi}{c^2}e^FR^2PV
\frac{\partial}{\partial m}, \\
\frac{\partial}{\partial r}&=4\pi\rho R^2R'\frac{\partial}{\partial m}.\label{drdm}
\end{align}
Here $(\partial/\partial t)_r$ stands for the partial derivative with respect to $t$ keeping $r$ constant, and $(\partial/\partial t)_m$ stands for that keeping $m$ constant.

Note that
\begin{equation}
\frac{\partial R}{\partial m}=\frac{1}{4\pi \rho R^2},
\end{equation}
and
\begin{equation}
\rho=\Big(4\pi R^2\frac{\partial R}{\partial m}\Big)^{-1}. \label{rhoR}
\end{equation}
Thus (\ref{Eq1}) reads
\begin{equation}
e^{-F}\Big(\frac{\partial R}{\partial t}\Big)_m=\Big(1+
\frac{P}{c^2\rho}\Big)V,\label{Ec1}
\end{equation}
and the equation (\ref{Eq4}) reads 
\begin{align}
e^{-F}\Big(\frac{\partial V}{\partial t}\Big)_m&=
\frac{4\pi}{c^2}R^2PV\frac{\partial V}{\partial m}-
GR\Big(\frac{m}{R^3}+\frac{4\pi P}{c^2}\Big) + \nonumber \\
&-\Big(1+\frac{V^2}{c^2}-
\frac{2Gm}{c^2R}\Big)
\Big(1+\frac{P}{c^2\rho}\Big)^{-1}\cdot
4\pi R^2\frac{\partial P}{\partial m}, \label{Ec2}
\end{align} 
where we have used the relation
$$\frac{P'}{R'}=4\pi\rho R^2\frac{\partial P}{\partial m},
$$
which comes from (\ref{drdm}).

Summing up the system of equations (\ref{Ec1})(\ref{Ec2}) should be solved, while $\rho, P=P(\rho)$ are given function of $R^2\partial R/\partial m$ through (\ref{rhoR}). Moreover under the assumption (A1) specified in the next section, we can put 
\begin{equation}
F=-\frac{u}{c^2}+F(0) \label{Fu}
\end{equation}
in order to fix the idea, where $F(0)$ is a constant and
$$u=\int_0^{\rho}
\frac{1}{\rho+P(\rho)/c^2}\frac{dP}{d\rho}d\rho$$
is a given function of $R^2\partial R/\partial m$, too. See (\ref{25}).
Hence the unknown functions are only
$(t, m)\mapsto R$ and $(t,m)\mapsto V$. 

The system of equations (\ref{Ec1})(\ref{Ec2}) will be called ($E_c$):
\begin{eqnarray*}
 e^{-F}\frac{\partial R}{\partial t}&=&\Big(1+
\frac{P}{c^2\rho}\Big)V,  \\
e^{-F}\frac{\partial V}{\partial t}&=&
\frac{4\pi}{c^2}R^2PV\frac{\partial V}{\partial m}-
GR\Big(\frac{m}{R^3}+\frac{4\pi P}{c^2}\Big) + \\
&-&\Big(1+\frac{V^2}{c^2}-
\frac{2Gm}{c^2R}\Big)
\Big(1+\frac{P}{c^2\rho}\Big)^{-1}\cdot
4\pi R^2\frac{\partial P}{\partial m}.  
\end{eqnarray*} 
Here we have written
$\displaystyle \frac{\partial R}{\partial t}, 
\frac{\partial V}{\partial t}$ simply instead of
$\displaystyle \Big(\frac{\partial R}{\partial t}\Big)_m, 
\Big(\frac{\partial V}{\partial t}\Big)_m$. The non-relativistic limit as $c\rightarrow +\infty$ is of course ($E_{\infty}$):
\begin{align*}
 \frac{\partial R}{\partial t}&=V, \\
\frac{\partial V}{\partial t}&=-\frac{Gm}{R^2}-4\pi
R^2\frac{\partial P}{\partial m},
\end{align*}
which is reduced to the second-order single equation 
\cite[\!(4)\!]{MakinoOJM}, where $g_0, r$ stand for $G, R$.\\

{\bf Supplementary Remark 1}.
The function $F=F(t,r)$ in the components of the metric (2.5) 
should satisfy (2.9).
Therefore, generally speaking, the formula (2.27) should read
$$
F=-\frac{u}{c^2}+F_+(t),
$$
where $F_+(t)$ is an arbitrary smooth function of $t$, being
constant with respect to $r$, or,
$$
e^{2F}=C(t)^2\kappa e^{-2u/c^2},
$$
where
$\kappa$ is a positive constant which will be specified in the next section, (see (3.9)),
and $C(t)$ is an arbitrary positive smooth function of $t$.
Then the left-hand sides of (2.25), (2.26) or ($E_c$) should be interpreted with
$$
e^{-F}\frac{\partial}{\partial t}=
\frac{1}{C(t)}\frac{1}{\sqrt{\kappa}}e^{u/c^2}\frac{\partial}{\partial t}.
$$
Of course we can and shall assume that $C(t)\equiv 1$ by taking
$$
t^*=t^*(t):=\int_0^tC(t')dt'
$$
instead of $t$, that is, we specify
\begin{equation}
e^F=\sqrt{\kappa}\exp \Big(-\frac{u}{c^2}\Big), \label{C=1}
\end{equation} 
without loss of generality.

In this sense, 
if we are allowed to forestall the discussion, we should say that,
in order to fix the idea, the definitions of
$J, H_1, H_2$ in the following Section 6 (see (6.8a), (6.8b) and (6.9)) should be done
by using (\ref{C=1}),
where $u$ is a given function of $\rho$ given by
$$\rho=\bar{\rho}(1+y)^{-2}
\Big(1+y+r\frac{\partial y}{\partial r}\Big)^{-1}. \eqno(4.5)
$$

\section{Equilibrium configurations}

Let us consider a solution of (\ref{Eq1})-(\ref{Eq4}) which is independent of $t$, that is, $F=F(\rho(r)), H=H(r), \rho=\rho(r), P=P(\rho(r)), V\equiv 0, R\equiv r$. Then the system of equations (\ref{Eq1})-(\ref{Eq4}) are reduced to 
$$0=Gr\Big(\frac{m}{r^3}+\frac{4\pi P}{c^2}\Big)+
\Big(1-\frac{2Gm}{c^2r}\Big)
\frac{P'}{\rho+P/c^2}.
$$
Therefore the equation to be studied is
\begin{subequations}
\begin{eqnarray}
\frac{dm}{dr}&=&4\pi r^2\rho, \label{TOV1} \\
\frac{dP}{dr}&=&-(\rho+P/c^2)
\frac{G(m+4\pi r^3P/c^2)}{r^2(1-2Gm/c^2r)}. \label{TOV2}
\end{eqnarray}
\end{subequations}
This equation was first derived by Oppenheimer-Volkoff \cite{OppenheimerV} in 1939.\\

Let us observe solutions of the Tolman-Oppenheimer-Volkoff equation (3.1). We assume (A0). 

\begin{Proposition}
Let $\rho_c(>0)$ and $P_c=P(\rho_c)$ be given central density and central pressure. Then there is a unique local solution $(m(r), P(r)), 0\leq r\leq \delta, of (3.1), \delta$ being a small positive number, such that $m=0, P=P_c$ at $r=0$. Moreover we have
\begin{align*}
m&=\frac{4\pi}{3}\rho_cr^3+O(r^5), \\
P&=P_c-(\rho_c+P_c/c^2)
G(4\pi\rho_c/3+4\pi P_c/c^2)\frac{r^2}{2}+O(r^4)
\end{align*}
as $r\rightarrow 0$.
\end{Proposition}

A proof can be found in \cite{Makino98}.

We consider the domain of the equation (3.1) as $\mathcal{D}:=\{
(r,m,P)\ |\   0<r<+\infty, 0<P<+\infty,
0<2Gm/c^2r<1\}$. Prolonging the local solution as long as possible in the domain 
$\mathcal{D}$, we have $(0,r_+)$ the maximal interval of existence. Here $r_+\leq+\infty$ is a constant. 

\begin{Definition}
If $r_+=+\infty$, the solution will be called a long equilibrium with the central density $\rho_c$. If $r_+<+\infty$, the solution will be called a short equilibrium.
\end{Definition}

{\bf Remark.}    It will be shown that, if $r_+<+\infty $,
$\rho$ and $P$ tend to $0$ but $2Gm/c^2r$
tends to a positive number strictly less than $1$
as $r\rightarrow r_+-0$. In this sense the solution can be said 
to be `short' if $r_+<+\infty$.\\

The equation of state for neutron stars is given by
\begin{align*}
P&=Kc^5\int_0^{\zeta}\frac{q^4dq}{(1+q^2)^{1/2}} \\
&=\frac{3}{8}Kc^5\Big(
5(1+\zeta^2)(\frac{2}{3}\zeta^2-1)+\log(\zeta+(\zeta^2+1)^{1/2})\Big) \\
\rho&=3Kc^3\int_0^{\zeta}(1+q^2)^{1/2}q^2dq \\
&=\frac{3}{8}Kc^3((2\zeta^2+1)(\zeta(\zeta^2+1)^{1/2}-
\log(\zeta+(\zeta^2+1)^{1/2})).
\end{align*}
See \cite[p. 188, (6.8.4), (6.8.5)\!]{ZeldovichN}. In this case we have
$$P=\frac{1}{5}K^{-2/3}\rho^{5/3}(1+[K^{-2/3}\rho^{2/3}/c^2]_1),$$
where $[X]_1$ stands for a convergent power series of the form
$\sum_{j\geq 1}a_jX^j$.Keeping in mind this case, we suppose the following assumption of the behavior of $P(\rho)$ ar $\rho \rightarrow 0$:\\

(A1)\hspace{10mm} {\bf There are positive constants $A, \gamma$ such that
$$P=A\rho^{\gamma}(1+[\rho^{\gamma-1}]_1)$$
as $\rho\rightarrow +0$, and $1<\gamma <2$. }\\

Under the assumptions (A0)(A1) we can introduce the new variable $u$ by
\begin{equation}
u=\int_0^P\frac{dP}{\rho+P/c^2},
\end{equation}
which satisfies
$$u=\frac{A\gamma}{\gamma-1}\rho^{\gamma-1}(1+[\rho^{\gamma-1}]_1)$$
as $\rho\rightarrow +0$. Let $(m(r), P(r)), 0<r<r_+,$ be an equilibrium, where $(0,r_+)$ is the maximal interval of existence. Then the corresponding $u=u(r)$ satisfies
\begin{equation}
r\frac{du}{dr}=-\frac{G(m+4\pi r^3P/c^2)}{r(1-2Gm/c^2r)}.
\end{equation}
Then $u(r)$ is monotone decreasing, and moreover we have

\begin{Proposition}
$u(r)\rightarrow 0$ as $r\rightarrow r_+-0$.
\end{Proposition}

Proof is the same as that of \cite[Lemma]{Makino98}. (We do not use the assumption $\gamma>4/3$.)\\

Let us introduce the variables
\begin{equation}
x=\frac{m}{ur},\qquad y=
4\pi r^2\frac{\rho^2}{P}.
\end{equation}
The equations read
\begin{align}
r\frac{dx}{dr}&=\alpha(u)-x+x^2\tilde{G}, \\
r\frac{dy}{dr}&=y(2-\beta(u)x\tilde{G}), \\
r\frac{du}{dr}&=-ux\tilde{G},
\end{align}
where
\begin{align*}
\alpha&=\frac{P}{u\rho}=\frac{\gamma-1}{\gamma}+[u]_1, \\
\beta&=\Big(2\frac{dP}{d\rho}-\frac{u}{P}\Big)
=\frac{2-\gamma}{\gamma-1}+[u]_1, \\
\tilde{G}&=\frac{G(1+4\pi r^3P/mc^2)}{1-2Gm/rc^2}
=\frac{G(1+\omega(u)y/c^2x)}{1-2Gux/c^2}, \\
\omega&=\frac{P^2}{u\rho^2}=[u]_1.
\end{align*}

\begin{Proposition}
Let $x(r)$ be corresponding to an equilibrium $(m(r), P(r)), 0<r<r_+$. If there is $r_0\in (0, r_+)$ such that $x(r_0)>1/G$, then $r_+<+\infty$ and enjoys the estimate
$$r_+<r_0\exp\Big(\frac{1}{Gx(r_0)-1}\Big).$$
\end{Proposition}

A proof can be found in the last part of the proof of \cite[Theorem 1]{Makino98}.\\

As in \cite{Makino98} we can claim
\begin{Proposition}
If $4/3<\gamma<2$, then any equilibrium is short.
\end{Proposition}

When $6/5<\gamma \leq 4/3$, it is known that, if $A$ is small and if $P(\rho)$ is sufficiently near to the exact $\gamma$-law $P=A\rho^{\gamma}$, then any equilibrium is short. See \cite{RendallS}. 
Even if $1<\gamma \leq 6/5$, it is possible that there are short equilibria, since Proposition 3 guarantees existence of tails of short equilibria in any case and we can arbitrarily modify the equation of state in the higher density region. Anyway in this article we assume (A0)(A1) only with $1<\gamma<2$ and suppose that a short equilibrium is given in front of us.\\

Let us observe roughly the behavior of a short equilibrium $(m(r), P(r))$ at the surface $r=r_+$.

By Proposition 2 we have $u \in C((0, r_+])$ with $u(r_+)=0$ and $P(r), \rho(r)$ are so, too. Hence
$$r\mapsto m(r)=\int_0^r4\pi r'^2\rho(r')dr'$$
belongs to $C((0, r_+])$. Put 
\begin{equation}
m_+=m(r_+)=\int_0^{r_+}4\pi r^2\rho(r)dr. \label{M}
\end{equation}
By definition we have $1-2Gm/c^2r>0$. Therefore
\begin{equation}
\kappa=\lim_{r\rightarrow r_+}
1-2Gm/c^2r=
1-2Gm_+/c^2r_+ \label{kappa}
\end{equation}
is non-negative. We claim that $\kappa >0$. Otherwise, if $\kappa=0$,
$$\frac{d}{dr}(1-2Gm/c^2r)
=-\frac{2G}{c^2}\Big(4\pi r\rho-\frac{m}{r^2}\Big)\rightarrow \frac{2Gm_+}{c^2r_+^2}=
\frac{1}{r_+}$$
as $r\rightarrow r_+-0$ and
$$1-2Gm/c^2r \sim -\frac{1}{r_+}(r_+-r),$$
which contradicts to $1-2Gm/c^2r>0$ for $r<r_+$. Hence $\kappa>0$ and
$$\frac{du}{dr}\rightarrow -K$$
as $r\rightarrow r_+-0$. Here
\begin{equation}
K=\frac{Gm_+}{r_+^2(1-2Gm_+/c^2r_+)} \label{K}
\end{equation}
is a positive constant. Hence, since $u\rightarrow 0$ as $r\rightarrow r_+$, we see
$$u\sim K(r_+-r) $$
and thus we have

\begin{Proposition}
Let $ (m(r), P(r)), 0<r<r_+,$ be a short equilibrium. Then we have
$$\rho\sim \Big(\frac{(\gamma-1)K}{A\gamma}\Big)^{\frac{1}{\gamma-1}}
(r_+-r)^{\frac{1}{\gamma-1}}$$
as $r\rightarrow r_+-0$, where $K$ is the positive constant given by (\ref{K}).
\end{Proposition}

{\bf Remark.}  If $(m(r),P(r)), 0<r<r_+,$ is a short equilibrium, then, on $r\geq r_+$, we put $\rho=P=0$  (vacuum), and
$$ds^2=\Big(1-\frac{2Gm_+}{c^2r}\Big)c^2dt^2
-\frac{dr^2}{1-\frac{2Gm_+}{c^2r}}-r^2(d\theta^2
+\sin^2\theta d\phi^2),$$
which is the Schwarzschild's metric. See \cite[p.301]{LandauL}.
Here we must take 
$$F(0)=\frac{1}{2}\log \kappa=\frac{1}{2}\log\Big(1-
\frac{2Gm_+}{c^2r_+}\Big).$$
Then the components of the metric are continuously differentiable across $r=r_+$.\\

More precise behavior of the equilibrium at the surface can be given as follows.

\begin{Proposition}
Assume (A0)(A1), and let $(m(r), P(r)), 0<r<r_+,$ be a short equilibrium.
If $\displaystyle\frac{\gamma}{\gamma-1}$ is an integer, then
$u(r)$
is analytic at $r=r_+$. 
\end{Proposition}

Proof. We consider the variables
$$X=\frac{1}{x}=\frac{ur}{m},\qquad
Y=\frac{y}{x^2}=\frac{4\pi r^4u^2\rho^2}{m^2P}.$$
Since $du/dr<0$, we can take $u$ as the independent variable instead of $r$, and the equations turn out to be
\begin{subequations}
\begin{align}
u\frac{dX}{du}&=\Big(1+\frac{1}{\tilde{G}}(-X+\alpha Y)\Big)X \label{eqX} \\
u\frac{dY}{du}&=
\Big(2+\beta+\frac{1}{\tilde{G}}(-4X+2\alpha Y)\Big)Y, \label{eqY}
\end{align}
\end{subequations}
where we note
\begin{equation}
\tilde{G}=G\Big(1+\frac{\omega u}{c^2}
\frac{Y}{X}\Big)\Big/
\Big(1-\frac{2G}{c^2}\frac{u}{X}\Big)
\end{equation}

Note that $\tilde{G}>0$ and $\tilde{G}\rightarrow G/\kappa$ as $u\rightarrow +0$,
where $\kappa$ is the positive constant given in (\ref{kappa}). Put
\begin{equation}
\check{X}=\frac{X}{u}=\frac{r}{m},\qquad
\check{Y}=\frac{Y}{u^{\frac{\gamma}{\gamma-1}}}=
\frac{4\pi r^4u^{\frac{\gamma-2}{\gamma-1}}\rho^2}{m^2P}.
\end{equation}
We know that $u\mapsto \check{X}$ and $u\mapsto \check{Y}$ belong to
$C([0, u_c))$ and $\check{X}|_{u=0}, \check{Y}|_{u=0}$ are positive.
Therefore $u\mapsto \tilde{G}=\displaystyle G\Big(1+
\frac{\omega u^{1/(\gamma-1)}}{c^2}
\frac{\check{Y}}{\check{X}}\Big)\Big/
\Big(1-\frac{2G}{c^2}\frac{1}{\check{X}}\Big)$
belongs to $C([0, u_c))$.

Integrating (\ref{eqX}), we have
$$X=C_1u\exp\Big[\int_0^u\frac{1}{\tilde{G}}(-\check{X}+
\alpha u^{\frac{\gamma}{\gamma-1}}\check{Y})du\Big].$$
Since the integrand is continuous, we see $u\mapsto
\check{X}$ belongs to $C^1([0, u_c))$. Integrating (\ref{eqY}), we have
$$Y=C_2u^{\frac{\gamma}{\gamma-1}}
\exp\Big[\int_0^u\Big(
\frac{1}{\tilde{G}}(-4\check{X}+2\alpha u^{\frac{\gamma}{\gamma-1}}\check{Y})+\Omega(u)\Big)du\Big],
$$
where
$$2+\beta=\frac{\gamma}{\gamma-1}+\Omega(u)u,\qquad \Omega(u)=[u]_0.$$

Fixing $u_0>0$ small, we put $\check{X}_0:=\check{X}(u_0),
\check{Y}_0:=\check{Y}(u_0)$. Since we know that
$$\check{X}(u)\rightarrow \check{X}_*:=\frac{r_+}{m_+},\qquad
\check{Y}(u)\rightarrow
\check{Y}_*:=4\pi\Big(\frac{\gamma-1}{A\gamma}
\Big)^{\frac{2-\gamma}{\gamma-1}}\frac{r_+^4}{m_+^2}
$$
as $u\rightarrow 0$, we see that, if $u_0$ is sufficiently small, then
$\check{X}_0, \check{Y}_0$ is arbitrarily near to $\check{X}_*, \check{Y}_*$. Now $(\check{X}(u), \check{Y}(u))$ is the unique solution of the integral equation
\begin{subequations}
\begin{align}
\check{X}(u)&=
\check{X}_0\exp\Big[
-\int_u^{u_0}\frac{1}{\tilde{G}}(-\check{X}+\alpha u^{\frac{\gamma}{\gamma-1}}\check{Y})du\Big], \label{EqX0} \\
\check{Y}(u)&=\check{Y}_0
\exp\Big[-\int_u^{u_0}\Big(\frac{1}{\tilde{G}}
(-4\check{X}+2\alpha u^{\frac{\gamma}{\gamma-1}}\check{Y})
+\Omega(u)\Big)du\Big].
\label{EqY0}
\end{align}
\end{subequations}
Let us denote by $D_{\delta}$ the set 
$\{(z_1,z_2)\in \mathbb{C}^2\  |\  |z_1-\check{X}_*|<\delta,
|z_2-\check{Y}_*|<\delta\}$, $\delta$ being small positive number. Note that, if
$|u|\leq \varepsilon_0, \varepsilon_0$ being a fixed small positive number, and if $(\check{X}, \check{Y}) \in D_{\delta}, 0<\delta\leq \delta_0$, then
$|1/\tilde{G}|\leq M_0, M_0$ depending upon only $\varepsilon_0, \delta_0$.
In fact, since $\check{X}_*>0$ and $\delta$ is very small, we can suppose that $(\check{X},\check{Y}) \in D_{\delta}$ guarantees $|\check{X}|\geq \delta$. Let us consider the functional family $\mathfrak{F}(\varepsilon, \delta)$ of all analytic functions $(\phi_1(u), \phi_2(u))$ defined and analytic
for $|u|\leq \varepsilon$ and valued in $D_{\delta}$. The right-hand side of (\ref{EqX0}), (\ref{EqY0}), in which $\check{X}=\phi_1(u), \check{Y}=\phi_2(u)$, will be denoted by $\bar{\phi}_1(u), \bar{\phi}_2(u)$. Then it is easy to see that,
if $|u_0|\leq\varepsilon, \varepsilon$ being sufficiently small, then
$(\check{X}_0, \check{Y}_0) \in D_{\delta/2}$, and if $(\phi_1, \phi_2) \in
\mathfrak{F}(\varepsilon, \delta)$, then $(\bar{\phi}_1, \bar{\phi}_2) \in
\mathfrak{F}(\varepsilon, \delta)$. Applying the well-known fixed point theorem
(see, e.g. \cite[Chapter I, Th\'{e}or\`{e}me 7]{HukuharaKM}), we have a fixed function in $\mathfrak{F}(\varepsilon, \delta)$. This is our
$(\check{X}(u), \check{Y}(u))$ by dint of the uniqueness.

 Integrating
$$\frac{1}{r}\frac{dr}{du}=-\frac{\check{X}}{\tilde{G}},$$
we see $u\mapsto r$ is analytic and $dr/du <0$ including $u=0$. Hence the inverse function $r\mapsto u$ is analytic at $r=r_+$. $\square$\\

Hereafter we suppose\\ 

(A2) \hspace{10mm}     {\bf $1<\gamma<2$ and $\displaystyle \frac{\gamma}{\gamma-1}$ is an integer.}\\

Under this assumption (A2), $\frac{1}{\gamma-1}$ is an integer, and,
since
\begin{equation}
\rho=\Big(\frac{\gamma-1}{A\gamma}u\Big)^{\frac{1}{\gamma-1}}(1+
[u]_1), \label{rho.u}
\end{equation}
the density distribution $\rho$ of the equilibrium is analytic at $r=r_+$, too:
\begin{equation}
\rho=\Big(\frac{(\gamma-1)K}{A\gamma}\Big)^{\frac{1}{\gamma-1}}
(r_+-r)^{\frac{1}{\gamma-1}}
(1+[r_+-r]_1).
\label{rho.r+}
\end{equation}

{\bf Supplementary Remark 2}.
Here we are considering a short
equilibrium with surface $r=r_+$ and the Schwartzschild's metric in the
exterior vacuum region. In other words, the metric 
$$ds^2=g_{00}c^2dt^2+
g_{11}dr^2-r^2(d\theta^2+
\sin^2\theta d\phi^2)$$
is given by
\begin{align*}
g_{00}&= \begin{cases}
e^{2F}=\kappa e^{-2u/c^2} \qquad (0\leq r \leq r_+) \\
\displaystyle 1-\frac{2Gm_+}{c^2r} \qquad (r_+<r), 
\end{cases}\\
-g_{11}&=\begin{cases}
\displaystyle  e^{2H}=\Big(1-\frac{2Gm}{c^2r}\Big)^{-1} 
\qquad (0\leq r\leq r_+) \\
\\
\displaystyle \Big(1-\frac{2Gm_+}{c^2r}\Big)^{-1} \qquad (r_+<r).
\end{cases}
\end{align*}
Let us see that the components $g_{00}, g_{11}$ are of class $C^2$
across $r=r_+$.

It is clear that
$g_{00}$ and $g_{11}$ are continuous since $u \rightarrow 0, m\rightarrow
m_+$ as $r\rightarrow r_+-0$ and $\kappa=1-2Gm_+/c^2r_+$.
Moreover we have
$$
\frac{d}{dr}
g_{00}\Big|_{r=r_+-0}=
-\frac{2\kappa}{c^2}\frac{du}{dr}\Big|_{r=r_+-0}=
\frac{2Gm_+}{c^2r_+^2}, $$
since
$du/dr \rightarrow -K$ with $K$ given by (3.10),
and
$$
\frac{d^2}{dx^2}
g_{00}\Big|_{r=r_+-0}=
\frac{4\kappa}{c^4}\Big(\frac{du}{dr}\Big)^2_{r=r_+-0}
-\frac{2\kappa}{c^2}\frac{d^2u}{dr^2}\Big|_{r=r_+-0}
=-\frac{4Gm_+}{c^2r_+^3}.
$$
This can be verified by differentiating the equation (3.3) and
seeing
$$\frac{d^2u}{dr^2}\Big|_{r=r_+-0}=
\frac{2Gm_+}{r_+^3\kappa}+
\frac{2}{c^2}\Big(\frac{Gm_+}{r_+^2\kappa}\Big)^2. $$
Hence $g_{00}$ is twice continuously
differentiable at $r=r_+$.
On the other hand, it is easy to see that the patched function
$$\tilde{m}(r)=\begin{cases}
m(r) & \qquad (0\leq r\leq r_+) \\
m_+& \qquad (r_+<r)
\end{cases}$$
is of class $C^k$ iff $\gamma<k/(k-1)$ since
$$\frac{dm}{dr}=4\pi\rho r^2=
4\pi r_+^2\Big(\frac{(\gamma-1)K}{A\gamma}\Big)^{1/(\gamma-1)}
(r_+-r)^{\frac{1}{\gamma-1}}(1+[r_+-r]_1).$$
Hence $\tilde{m}$ and $g_{11}$ are of class $C^2$ since $\gamma<2$.

\section{Equations for perturbations}

Let us fix a short equilibrium $\rho(r)$ which is positive on $0\leq r<r_+$. Put $m_+=m(r_+)$. Then we can take $m$ as an independent variable and get an equilibrium $\rho=\bar{\rho}(m)$ and $r=r(m), 0\leq m\leq m_+$. We have to consider solutions of ($E_c$) near to this equilibrium of the form
\begin{align}
R&=r(m)(1+y), \label{Ry} \\
V&=r(m)v. \label{Vv}
\end{align}
Here $y$ and $v$ are small perturbations. The equations turn out to be
\begin{align}
e^{-F}\frac{\partial y}{\partial t}&=\Big(1+\frac{P}{c^2\rho}\Big)v, \label{pert1} \\
e^{-F}\frac{\partial v}{\partial t}&=
\frac{4\pi}{c^2}r^2(1+y)^2Pv\frac{\partial }{\partial m}(rv)+ \nonumber \\
&-\frac{G}{r^3(1+y)^2}\Big(m+
\frac{4\pi}{c^2} P r^3(1+y)^3\Big)+ \nonumber \\
&-\Big(1+\frac{r^2v^2}{c^2}-\frac{2Gm}{c^2r(1+y)}\Big)\Big(1+\frac{P}{c^2\rho}\Big)^{-1}\cdot
4\pi r(1+y)^2\frac{\partial P}{\partial m}. \label{pert2}
\end{align}
Instead of $m$, let us take $r=r(m)$ as the independent variable. Since
$$\frac{dm}{dr}=4\pi\bar{\rho}r^2,$$
we see
$$\frac{\partial}{\partial m}=\frac{1}{4\pi\bar{\rho}r^2}\frac{\partial}{\partial r}.$$
Therefore (\ref{rhoR}) and (\ref{Ry}) imply
\begin{equation}
\rho=\bar{\rho}(1+y)^{-2}\Big(1+
y+r\frac{\partial y}{\partial r}\Big)^{-1} \label{rho}
\end{equation}
so that
$$\rho=\bar{\rho}\Big(1-3y-r\frac{\partial y}{\partial r}
+\Big[y, r\frac{\partial y}{\partial r}\Big]_2\Big).
$$
Here $[X_1,X_2]_2$ denotes a convergent double power series of the form $$\sum_{k_1+k_2\geq 2}a_{k_1k_2}X_1^{k_1}X_2^{k_2}.$$

Let us recall that $ (\bar{\rho})^{\gamma-1}\in C^{\infty}([0, r_+])$, provided that $\gamma/(\gamma-1)$ is an integer, say, (A2). 

The equation (\ref{pert2}) reads
\begin{align}
e^{-F}\frac{\partial v}{\partial t}&=
\frac{1}{c^2}(1+y)^2\frac{P}{\bar{\rho}}v\frac{\partial}{\partial r}(rv)+ \nonumber \\ 
&-\frac{G}{r^3(1+y)^2}\Big(m+\frac{4\pi}{c^2}Pr^3(1+y)^3\Big)+ \nonumber \\
&-\Big(1+\frac{r^2v^2}{c^2}-\frac{2Gm}{c^2r(1+y)}\Big)
\Big(1+\frac{P}{c^2\rho}\Big)^{-1}
\frac{(1+y)^2}{r\bar{\rho}}\frac{\partial P}{\partial r}. \label{pert2b}
\end{align}
We have to solve (\ref{pert1})(\ref{pert2b}) for unknown functions $(t,r)
\mapsto y, v$, where $r$ is confined to the fixed interval $[0,r_+]$. Here $m=m(r)$ is determined by the equilibrium through
$$m=4\pi \int_0^r\bar{\rho}(r)r^2dr,
$$
and
$\rho, P(\rho), u(\rho)$ are given functions of
$\bar{\rho}(r)$ and the unknowns $y, r\partial y/\partial r$ through (\ref{rho}). \\

The perturbation of $\rho$ is expressed by (\ref{rho}). Similar expressions of $P$ and $u$ are necessary. If $P(\rho)$ was the exact $\gamma$-law, say, if
$P=A\rho^{\gamma}$, then we would have
\begin{align*}
P&=\bar{P}(1+y)^{-2\gamma}
\Big(1+y+r\frac{\partial y}{\partial r}\Big)^{-\gamma} \\
&=\bar{P}\Big(1-\gamma\Big(3y+r\frac{\partial y}{\partial r}\Big)+
\Big[y,r\frac{\partial y}{\partial r}\Big]_2\Big).
\end{align*}
However 
this exact $\gamma$-law is not treated by this article, since it violates the condition $dP/d\rho <c^2$ for large $\rho$. 
Our case should be treated more carefully.

We should introduce the quantity
\begin{equation}
\gamma^P:=\frac{\rho}{P}\frac{dP}{d\rho}.
\end{equation}
Then under the assumption (A1) we see
$$\gamma^P=\gamma +[u]_1$$
and, using this function, we can express
\begin{equation}
P=\bar{P}
\Big(1-\gamma^P(\bar{u})\Big(3y+r\frac{\partial y}{\partial r}\Big)
-\Phi^P\Big(\bar{u}, y, r\frac{\partial y}{\partial r}\Big)\Big),
\end{equation}
where $$\Phi^P(u, y, ry')=[u; y, ry']_{0;2}.$$
Here $[X_0;X_1,X_2]_{0;2}$ denotes a convergent triple power series of the form
$$\sum_{k_0\geq 0, k_1+k_2\geq 2}a_{k_0k_1k_2}
X_0^{k_0}X_1^{k_1}X_2^{k_2}.$$



We note that
\begin{align}
\Big(1+\frac{P}{c^2\rho}\Big)^{-1}&=
\Big(1+\frac{\bar{P}}{c^2\bar{\rho}}\Big)^{-1}
\Big(1+\frac{\bar{P}}{c^2\bar{\rho}}
\Big(1+\frac{\bar{P}}{c^2\bar{\rho}}\Big)^{-1}
(\gamma^P-1)\Big(3y+r\frac{\partial y}{\partial r}\Big)+ \nonumber \\
&+\Big[\!(\bar{u};y,r\frac{\partial y}{\partial r}\Big]_{0;2}\Big).
\end{align}

{\bf Supplementary Remark 3}.
The letter `$r$' is used, on the one hand, 
as one of the co-moving coordinates for the metric (2.5),
and, on the other hand, as one of the independent variables of the equation
(4.6). However these two
quantities denoted by the same letter `$r$' do not coincide if we consider moving solutions.
Therefore, in order to clarify the relation between these two
quantities, we shall denote by $r^*$ the latter  $r$, say one of the independent variables for (4.6). 

In other words, 
the definition of $r^*=\varphi(t,r)$ is as following:
Put 
$$m=f_1(r):=4\pi\int_0^r\bar{\rho}(r')r'^2dr' \quad (0\leq r\leq r_+) $$
along the equilibrium fixed. Then we have
the inverse function $r=f_1^{-1}(m)$ defined on
$0\leq m\leq  m_+$. But along the moving solutions $m$ is one of the
variables of the equations (2.25)(2.26) or $(E_c)$ defined by
(2.11). So we denote
$$m=f_2(t,r):=
4\pi\int_0^r\rho(t,r')R(t,r')^2
\partial_rR(t,r')dr' $$
along the moving solutions under consideration. Then
we put
$$r^*=\varphi(t,r):=f_1^{-1}(f_2(t,r)).
$$

Let us determine the function $\varphi(t,x)$.The function
 $m=f_2(t,r)$ should satisfy (2.20d), that is,
\begin{equation}
e^{-F}\frac{\partial m}{\partial t}=-\frac{4\pi}{c^2}RPV.\label{Eq6}
\end{equation}
The left-hand side of (\ref{Eq6}) is
$$\frac{1}{\sqrt{\kappa}}e^{u/c^2}Df_1(r^*)\frac{\partial r^*}{\partial t}
=\frac{1}{\sqrt{\kappa}}
e^{u/c^2}\cdot 4\pi\bar{\rho}(r^*)
(r^*)^2\frac{\partial r^*}{\partial t}. $$
On the other hand, we are going to construct moving solutions of the
form
\begin{align*}
&R=r^*(1+y(t,r^*)), \qquad
V=r^*v(t, r^*), \\
&P=P(\rho), \qquad \rho=\bar{\rho}(r^*)
(1+y(t,r^*))^{-2}
\Big(1+y+r^*\frac{\partial y}{\partial r^*}\Big)^{-1}
\end{align*}
with $y, v \in C^{\infty}([0,T]\times [0, r_+])$,
which are very small. 
Suppose that we have constructed such solutions. Then the function
$\varphi(t,r)$ should satisfy
\begin{equation}
\frac{\partial}{\partial t}\varphi(t,r)=
-\frac{\sqrt{\kappa}}{c^2}e^{-u/c^2}
(1+y)^2\frac{P}{\bar{\rho}}v\cdot\varphi(t,r),\label{Eq7}
\end{equation}
where $u=u(\rho), y, P=P(\rho), \bar{\rho}, v$ in the right-hand
side  are evaluated at $(t,r^*)=(t, \varphi(t,r))$. The
formula
(\ref{Eq7}) can be considered as an ordinary differential
equation for $\varphi(\cdot, r)$ for each fixed $r$,
which determines $\varphi(\cdot, r)$ provided that
the initial value $\varphi(0, r)=f_1^{-1}(f_2(0,r))$ is given. 

But we can assume that $\varphi(0,r)=r$ without loss of generality.
In fact, $\varphi(0,r)=r$ means $f_1(r)=f_2(0,r)$ and,
even if $f_1\not= f_2(0,\cdot)$, we
can find the change of variable
$r=\psi(r^{\flat})$ such that
$f_1(r^{\flat})=f_2(0,\psi(r^{\flat}))$. Considering
$r^{\flat}$ instead of $r$, we can assume $f_1(r)=
f_2(0,r)$ or $\varphi(0,r)=r$.
Clearly the $C^{\infty}$-solution $\varphi$ is 
uniquely determined and $\varphi(t,r)-r$ is very
small with its derivatives. Of course $\varphi(t,0)=0$ and
$\varphi(t,r_+)=r_+$ since $P/\bar{\rho}$ vanishes at $r=r_+-0$.
Hence we have the solutions
$$
R=\varphi(t,r)(1+y(t,\varphi(t,r))), \qquad
V=\varphi(t,r)v(t,\varphi(t,r))
$$
and so on as functions of the original co-moving coordinates $t, r$. 

\section{Analysis of the linearized equation}

We are going to analyze the linearized equations to (\ref{pert1})(\ref{pert2b})
and establish the existence of time periodic 
solutions to the linearized equations of the form
$$y=\mbox{Const.}\sin(\sqrt{\lambda}t+
\mbox{Const.})\tilde{\psi}(r), $$
where $\lambda >0$ and $\tilde{\psi}(r)$ is an analytic function
 of $r$ in a neighborhood of $[0,r_+]$.\\

Using the formulas listed in the last part of the preceding section, we see that the linearizations of the equations (\ref{pert1})(\ref{pert2b}) turn out to be
\begin{subequations}
\begin{align}
e^{-F}\frac{\partial y}{\partial t}&=
\Big(1+\frac{P}{c^2\rho}\Big)v \label{linEq1} \\
e^{-F}\frac{\partial v}{\partial t}&=
E_2y''+
E_1y'+E_0y, \label{linEq2}
\end{align}
\end{subequations}
where $y''=\partial^2y/\partial r^2, y'=\partial y/\partial r$ and
\begin{subequations}
\begin{align}
E_2&=e^{-2H}(\rho+P/c^2)^{-1}\gamma^PP,\label{linEq22} \\
\frac{E_1}{E_2}&=
\frac{4\pi G}{c^2}e^{2H}(\rho+P/c^2)r-
(\rho+P/c^2)^{-1}\Big(1-\frac{1}{\gamma_p}\Big)\frac{P'}{c^2} 
+\frac{3}{r}+\frac{(\gamma^PPr)'}{\gamma^PPr}= \nonumber \\
&=F'+H'-
\frac{(1+P/c^2\rho)'}{1+P/c^2\rho}+
\frac{3}{r}+
\frac{(\gamma^PPr)'}{\gamma^PPr}, \label{linEq21} \\
E_0&=\frac{4\pi G}{c^2}\cdot 3(\gamma^P-1)P+ \nonumber \\
&+\Big(-1-3\gamma^Pe^{-2H}+
3(\gamma^P-1)e^{-2H}(1+P/c^2\rho)^{-1}\Big)
(\rho+P/c^2)^{-1}
\frac{P'}{r}+ \nonumber \\
&
+3e^{-2H}(\rho+P/c^2)^{-1}
\frac{(\gamma^PP)'}{r}
.  \label{linEq20}
\end{align}
\end{subequations}
Here $\rho, P, \gamma^P, F, H$ are abbreviations
for
the quantities
$\bar{\rho}(r), \bar{P}=P(\bar{\rho}(r)), 
\gamma^P(\bar{u}(r)), \\
\displaystyle \bar{F}=
F(\bar{u}(r))=-\frac{1}{c^2}\bar{u}(r)+\frac{1}{2}\log
\kappa, \bar{H}=-\frac{1}{2}\log\Big(1-\frac{2Gm}{c^2r}\Big)$ along the considered equilibrium. Throughout the above manipulations
we have used the equation
$$\frac{4\pi G}{c^2}e^{2H}(\rho+P/c^2)r=F'+H', $$
which can be derived from the differentiation of
$$e^{-2H}=1-\frac{2Gm}{c^2r}$$
and (\ref{TOV2}), and also the relation
$$\frac{(1+P/c^2\rho)'}{1+P/c^2\rho}=
\frac{1}{\rho+P/c^2}\Big(1-\frac{1}{\gamma^P}\Big)\frac{P'}{c^2}.$$

In other words, the linearized second order single equation is:
\begin{equation}
\frac{\partial^2y}{\partial t^2}+\mathcal{L}y=0, \label{linEq}
\end{equation}
where
\begin{equation}
\mathcal{L}y=
-\frac{a}{b}y''-\frac{a'}{b}y'+Qy=-\frac{1}{b}(ay')'+Qy,
\label{linDO}
\end{equation}
\begin{subequations}
\begin{align}
a&=\exp\Big[\int^r\frac{E_1}{E_2}dr\Big] =\frac{\gamma^PPr^4}{1+P/c^2\rho}e^{F+H}, \label{linDOa} \\
b&=(1+P/c^2\rho)^{-1}\rho r^{4}e^{-F+3H},\label{linDOb}\\
Q&=-e^{2F}(1+P/c^2\rho)E_0. \label{linDOc}
\end{align}
\end{subequations}

In order to investigate the spectral property of the second order linear differential operator $\mathcal{L}$, we reduce the eigenvalue problem
\begin{equation}
\mathcal{L}y=\lambda y
\end{equation}
to the normal form
\begin{equation}
-\frac{d^2\eta}{d\xi^2}+q(\xi)\eta=\lambda \eta
\end{equation}
by the Liouville transformation
\begin{subequations}
\begin{align}
\xi&=\int_0^r\sqrt{\frac{b}{a}}dr=
\int_0^r\sqrt{\frac{\rho}{\gamma^PP}}e^{-F+H}dr, \\
\eta&=(ab)^{1/4}y=
(\gamma^P\rho P)^{1/4}r^2
(1+P/c^2\rho)^{-1/2}e^Hy,
\end{align}
\end{subequations}
when the result is
\begin{equation}
q=Q+\frac{a}{4b}
\Big(\Big(\frac{a'}{a}+\frac{b'}{b}\Big)'-
\frac{1}{4}\Big(\frac{a'}{a}+\frac{b'}{b}\Big)^2+
\frac{a'}{a}\Big(\frac{a'}{a}+\frac{b'}{b}\Big)\Big).
\end{equation}
See \cite[p. 275, Theorem 6]{BirkhoffR}.

Since 
$$\sqrt{\frac{\rho}{\gamma^PP}}\sim \mbox{Const.}
(r_+-r)^{-1/2}, $$
we can define the finite value
\begin{equation}
\xi_+:=\int_0^{r_+}
\sqrt{\frac{\rho}{\gamma^PP}}e^{-F+H}dr.
\end{equation}
The interval $(0,r_+)$ is mapped onto $(0,\xi_+)$.

First let us observe the behavior of $q$ as $\xi \rightarrow 0 (r
\rightarrow 0)$. We see that $Q=O(1)$,
$$\frac{a'}{a}\sim \frac{4}{r},\quad \frac{a'}{a}+\frac{b'}{b}\sim \frac{8}{r},
\quad \Big(\frac{a'}{a}+\frac{b'}{b}\Big)'\sim
-\frac{8}{r^2}, $$
therefore
$$q\sim 2\gamma^PP\rho^{-1}e^{2F-2H}\Big|_{r=0}\frac{1}{r^2}.$$
On the other hand we have
$$\xi \sim \Big(\gamma^PP\rho^{-1}e^{2F-2H}\Big|_{r=0}\Big)^{-1/2}r.
$$
Hence we have
$$q \sim \frac{2}{\xi^2}.$$
Note that $2 > 3/4$. 

Next we observe the behavior of $q$ as $\xi\rightarrow \xi_+
(r\rightarrow r_+)$. Note that
$P'/\rho \rightarrow -K$, where $K$ is the constant defined by (\ref{K}). Therefore we see that
$Q=O(1)$. Moreover we have
$$\frac{\rho}{\gamma^P}\frac{d}{d\rho}\gamma^P=O(u) \rightarrow 0,$$
so that $(\gamma^P)'/\gamma^P=o(\rho'/\rho)$. Hence we see that
\begin{align*}
&\frac{a'}{a}\sim -\frac{\gamma}{\gamma-1}\frac{1}{r_+-r}, \quad 
\frac{a'}{a}+\frac{b'}{b}\sim -\frac{\gamma+1}{\gamma -1}\frac{1}{r_+-r}, \\
&\Big(\frac{a'}{a}+\frac{b'}{b}\Big)'\sim
-\frac{\gamma+1}{\gamma-1}\frac{1}{(r_+-r)^2},
\end{align*}
Therefore we have
$$q\sim
Ke^{2F-2H}\Big|_{r=r_+}\frac{(\gamma+1)(3-\gamma)}{16(\gamma-1)}
\frac{1}{r_+-r}.$$
On the other hand we have
$$\xi_+-\xi\sim
\frac{2}{\sqrt{(\gamma-1)K}}e^{-F+H}\Big|_{r=r_+}\sqrt{r_+-r}.
$$
Hence we have
$$q\sim \frac{(\gamma+1)(3-\gamma)}{4(\gamma-1)^2}
\frac{1}{(\xi_+-\xi)^2}.
$$
It follows from $1<\gamma<2$ that
$$\frac{(\gamma+1)(3-\gamma)}{4(\gamma-1)^2}>\frac{3}{4}.
$$
Therefore the both boundary points $\xi=0, \xi_+$ are of limit
point type, and \cite[p. 159, Theorem X.10]{ReedS} gives the following conclusion, which is the same as \cite[Proposition 1]{MakinoOJM}:
\begin{Proposition}
The operator $\mathfrak{T}_0, \mathcal{D}(\mathfrak{T}_0)=
C_0^{\infty}(0, \xi_+),
\mathfrak{T}_0\eta=-\eta_{\xi\xi}+q\eta$,
in $L^2(0,\xi_+)$ has the Friedrichs extension $\mathfrak{T}$, a self-adjoint
operator, whose spectrum consists of simple eigenvalues
$\lambda_1<\cdots <\lambda_n<\cdots \rightarrow +\infty$. In other words, the operator
$\mathfrak{S}_0, \mathcal{D}(\mathfrak{S}_0)=
C_0^{\infty}(0,r_+),
\mathfrak{S}_0y=\mathcal{L}y$ in $L^2((0,r_+), bdr)$ has the Friedrichs extension $\mathfrak{S}$, a self-adjoint operator with eigenvalues $(\lambda_n)_n$.
\end{Proposition}

In order to investigate the structure of the linear operator $\mathcal{L}$,
we introduce the new independent variable $x$ instead of $r$ defined by
\begin{equation}
x:=\frac{\tan^2\theta}{1+\tan^2\theta} \quad\mbox{with}\quad
\theta:=\frac{\pi\xi}{2\xi_+}
=\frac{\pi}{2\xi_+}
\int_0^r
\sqrt{\frac{\rho}{\gamma^PP}}e^{-F+H}dr.
\end{equation}
The interval $[0, r_+]$ of the variable $r$ is mapped onto
$[0,1]$ of $x$, and we have
\begin{subequations}
\begin{align}
\frac{d}{dr}&=\frac{\pi}{\xi_+}\sqrt{x(1-x)}\sqrt{\frac{b}{a}}\frac{d}{dx}, 
\label{d/dr} \\
\frac{d^2}{dr^2}&=\Big(\frac{\pi}{\xi_+}\Big)^2\frac{b}{a}\Big(
x(1-x)\frac{d^2}{dx^2}+
\Big(\frac{1-2x}{2}+
\frac{\xi_+}{\pi}\sqrt{x(1-x)}
\sqrt{\frac{a}{b}}\frac{1}{2}\frac{a}{b}\frac{d}{dr}\Big(\frac{b}{a}\Big)\Big)
\frac{d}{dx}\Big).\label{d2/dr2}
\end{align}
\end{subequations}

We note 
\begin{equation}
r\frac{d}{dr}=x[\!(x)\!]\frac{d}{dx},
\end{equation}
where and hereafter $[\!(x)\!]$ {\bf denotes an analytic function of $x$ in a neighborhood of the interval
$[0,1]$}. In fact, (\ref{d/dr})
implies the following observations: as $r\rightarrow 0 (x\rightarrow 0)$, we see
\begin{equation}
r=\frac{\xi_+}{\pi}C_0\sqrt{x}(1+[x]_1) \quad
\mbox{with}\quad
C_0=2\sqrt{\frac{\gamma^PP}{\rho}}
e^{F-H}\Big|_{r=0},
\end{equation}
and
$$\sqrt{\frac{b}{a}}=\frac{2\xi_+}{\pi}\frac{\sqrt{x}}{r}(1+[x]_1),
$$ so that
$$r\frac{d}{dr}=2x(1+[x]_1)\frac{d}{dx};
$$
as $r\rightarrow r_+ (x\rightarrow 1)$, we see
\begin{equation}
1-x=\Big(\frac{\pi}{\xi_+}\Big)^2C_1(r_+-r)(1+[r_+-r]_1)
\quad\mbox{with}\quad
C_1=\frac{1}{(\gamma-1)\kappa^2K},
\end{equation}
(see (\ref{rho.r+}) and note 
$e^{F-H}=\kappa+
[r_+-r]_1$ with $\kappa=1-2Gm_+/c^2r_+$) and
$$\sqrt{\frac{b}{a}}=\frac{\pi}{\xi_+}\frac{C_1}{\sqrt{1-x}}
(1+[1-x]_1),$$
so that
$$r\frac{d}{dr}=
\Big(\frac{\pi}{\xi_+}\Big)^2C_1r_+(1+[1-x]_1)\frac{d}{dx}.
$$

Now we can write
\begin{equation}
\Big(\frac{\xi_+}{\pi}\Big)^2\mathcal{L}y=
-x(1-x)\frac{d^2}{dx^2}-B\frac{dy}{dx}+\Big(\frac{\xi_+}{\pi}\Big)^2Qy,
\end{equation}
where
\begin{equation}
B=\frac{1-2x}{2}+
\frac{\xi_+}{\pi}
\sqrt{x(1-x)}\sqrt{\frac{a}{b}}
\Big(\frac{1}{2}\frac{a}{b}\frac{d}{dr}\Big(\frac{b}{a}\Big)+\frac{1}{a}
\frac{da}{dr}\Big).
\end{equation}

\textbullet\  As $r\rightarrow 0 (x\rightarrow 0)$, we see
$$B=\frac{5}{2}+[x]_1.$$
For
\begin{align*}
\frac{1}{2}\frac{a}{b}\frac{d}{dr}\Big(\frac{b}{a}\Big)+\frac{1}{a}
\frac{da}{dr}&=
\frac{1}{2}\frac{(\gamma^PP\rho)'}{\gamma^PP\rho}+2H'+
\frac{4}{r}-\frac{(1+P/c^2\rho)'}{1+P/c^2\rho} = \\
&=\frac{4}{r}(1+[r^2]_1)
\end{align*}
and
\begin{align*}
\sqrt{\frac{a}{b}}&=\sqrt{\frac{\gamma^PP}{\rho}}e^{F-H}= \\
&=\Big(\sqrt{\frac{\gamma^PP}{\rho}}e^{F-H}\Big)\Big|_{r=0}
(1+[r^2]_1)
=\frac{\pi}{2\xi_+}
\frac{r}{\sqrt{x}}(1+[x]_1).
\end{align*}
Clearly
$$Q=e^{-F}(1+P/c^2\rho)E_0=[r^2]_0=[x]_0.$$

\textbullet\  As $r\rightarrow r_+ (x\rightarrow 1)$, we see
$$B=-\frac{\gamma}{\gamma-1}+[1-x]_1.
$$
For
\begin{align*}
\frac{1}{2}\frac{a}{b}\frac{d}{dr}\Big(\frac{b}{a}\Big)+\frac{1}{a}
\frac{da}{dr} &=-\frac{1}{2}\frac{\gamma+1}{\gamma-1}\frac{1}{r_+-r}(1+
[r_+-r]_1) \\
&=-\frac{\gamma+1}{\gamma-1}\Big(\frac{\pi}{\xi_+}\Big)^2
\frac{C_1}{1-x}(1+[1-x]_1)
\end{align*}
and
$$\frac{a}{b}=\frac{1}{C_1}(r_+-r)(1+[r_+-r]_1)=
\frac{1}{C_1^2}\Big(\frac{\xi_+}{\pi}\Big)^2
(1-x)(1+[1-x]_1).
$$
Clearly
$$Q=[u]_0=[r_+-r]_0=[1-x]_0.
$$

Summing up, we have the following conclusion, which is the same as
\cite[Proposition 3]{MakinoOJM}:
\begin{Proposition}
We can write
\begin{equation}
\Big(\frac{\xi_+}{\pi}\Big)^2\mathcal{L}y=
-x(1-x)\frac{d^2y}{dx^2}
-\Big(\frac{5}{2}(1-x)-
\frac{N}{2}x\Big)\frac{dy}{dx}+L_1(x)\frac{dy}{dx}
+L_0(x)y,
\end{equation}
where $L_1(x)=x(1-x)[\!(x)\!], L_0(x)=[\!(x)\!]$. Here $N$ is the parameter
defined by
\begin{equation}
N=\frac{2\gamma}{\gamma-1}\quad\mbox{or}\quad
\gamma=\frac{N}{N-2}.
\end{equation}
\end{Proposition}

The assumption (A2) reads that $N$ is an even integer $>4$. As long as we are concerned with investigation of the analytic structure of the operator $\mathcal{L}$, we may assume that $\xi_+=\pi$ without loss of generality. \\

Anyway, we have the following

\begin{Proposition}
Let $\lambda=\lambda_n$ be a positive eigenvalue and let $\psi$ be an associated eigenfunction
which belongs to
$L^2([0,1]; x^{\frac{3}{2}}(1-x)^{\frac{N}{2}-1}dx)$. Then
\begin{equation}
Y_1=\sin(\sqrt{\lambda}t+\Theta_0)\psi(x) \label{Y1}
\end{equation}
is a time periodic solution of the linearized problem (\ref{linEq}). 
\end{Proposition}

Thanks to Proposition 8, we can claim the following proposition on the analytic property of the eigenfunction same as  \cite[Proposition 4]{MakinoOJM}:

\begin{Proposition}
We have
\begin{subequations}
\begin{align}
\psi(x)&=c_0(1+[x]_1) \quad \mbox{as}\quad x\rightarrow 0 \\
\psi(x)&=c_1(1+[1-x]_1) \quad\mbox{as}\quad x\rightarrow 1.
\end{align}
\end{subequations}
Here $c_0, c_1$ are non-zero constants. Other independent solutions of $\mathcal{L}y=\lambda y$ do not belong to $L^2([0,1];
(1-x)^{\frac{N}{2}-1}dx)$ as $ x\sim 1$.
\end{Proposition}

Therefore $\psi(x)=[\!(x)\!]$ and $Y_1$ is an analytic function of
$t \in \mathbb{C}$ and $x$ on a neighborhood of $[0,1]$ independent of $t$.\\

Hereafter we fix such a time periodic function $Y_1$.

\section{Rewriting of the equations
(\ref{pert1})(\ref{pert2b}) using
the linear operator $\mathcal{L}$}

Let us go back to the system of equations (\ref{pert1})(\ref{pert2b}). 
In order to rewrite these equations using the linear operator $\mathcal{L}$, we shall use the following observations.

We are considering the perturbed $P$ such that
\begin{equation}
P=\bar{P}(1-\gamma^P(\bar{u})(3y+z)-\Phi^P(\bar{u},y,z)),
\end{equation}
where $z=r\partial y/\partial r$.

Then we have
\begin{align}
-\frac{1}{r\bar{\rho}}\frac{\partial P}{\partial r}&=
-\frac{1}{r\bar{\rho}}\frac{d\bar{P}}{dr}+
\Big(1+\frac{1}{\gamma^P}\partial_z\Phi^P\Big)
\frac{1}{r\bar{\rho}}
\frac{\partial}{\partial r}
(\bar{P}\gamma^P(3y+z))+ \nonumber \\
&+\frac{\bar{P}}{r\bar{\rho}}\cdot[Q0]+
\frac{1}{r\bar{\rho}}\frac{d\bar{P}}{dr}\cdot[Q1]
\end{align}
where
\begin{subequations}
\begin{align}
[Q0]&:=2(\gamma^P+\partial_z\Phi^P)(1+y)^{-1}\frac{z^2}{r}+
\frac{d\bar{u}}{dr}
\Big(\partial_u\Phi^P-
\frac{1}{\gamma^P}\frac{d\gamma^P}{du}(3y+z)\partial_z\Phi^P\Big),
\\
[Q1]&:=\Phi^P-(3y+z)\partial_z\Phi^P.
\end{align}
\end{subequations}
Here we have used the relation
\begin{equation}
(\partial_y-3\partial_z)\Phi^P=
2(\gamma^P+\partial_z\Phi^P)(1+y)^{-1}z.
\end{equation}

Let us analyze
\begin{equation}
\mbox{the right-hand side of (\ref{pert2b})}=[R1]+[R2],
\end{equation}
where
\begin{subequations}
\begin{align}
[R1]&:=-G(1+y)\Big(\frac{m}{r^3(1+y)^3}+\frac{4\pi P}{c^2}\Big)+ \nonumber \\
&-\Big(1+\frac{r^2v^2}{c^2}-\frac{2Gm}{c^2r(1+y)}\Big)
(1+P/c^2\rho)^{-1}\frac{(1+y)^2}{r\bar{\rho}}\frac{\partial P}{\partial r}, \\
[R2]&:=\frac{1}{c^2}(1+y)^2\frac{P}{\bar{\rho}}(v^2+vw) \quad\mbox{with}\quad w=r\frac{\partial v}{\partial r}.
\end{align}
\end{subequations}

Let us put
$$[R1]=[R3]+[R4]+[R5]+[R6]+[R7],$$
where
\begin{align*}
[R3]&:=-G
(1+y)
\Big(\frac{m}{r^3(1+y)^3}+\frac{4\pi P}{c^2}\Big) \\
&=-\frac{Gm}{r^3}-\frac{4\pi G\bar{P}}{c^2}+[R3L]+[R3Q], \\
[R3L]&:=\frac{Gm}{r^3}\cdot
2y
+\frac{4\pi G}{c^2}\bar{P}\gamma^P(3y+z)
+\frac{4\pi G}{c^2}\bar{P}y
, \\
[R3Q]&:=-\frac{Gm}{r^3}\Big(
\frac{1}{(1+y)^2}-(1-2y)\Big)+
\frac{4\pi G}{c^2}\bar{P}\Phi^P
+\frac{4\pi G}{c^2}(P-\bar{P})y
; 
\end{align*}
\begin{align*}
[R4]&:=-\Big(
1+\frac{r^2v^2}{c^2}-
\frac{2Gm}{c^2r(1+y)}\Big)
(1+P/c^2\rho)^{-1}
\frac{(1+y)^2}{r\bar{\rho}}\frac{d\bar{P}}{dr} \\
&=-\Big(1-\frac{2Gm}{c^2r}\Big)(1+\bar{P}/c^2\bar{\rho})^{-1}
\frac{1}{r\bar{\rho}}\frac{d\bar{P}}{dr}+
[R4L]+[R4Q], \\
[R4L]&:=
\Big(-\frac{2Gm}{c^2r}(1+\bar{P}/c^2\bar{\rho})^{-1}\cdot y+
\Big(1-\frac{2Gm}{c^2r}\Big)\frac{\bar{P}}{c^2\bar{\rho}}
(1+\bar{P}/c^2\rho)^{-1}(\gamma^P-1)(3y+z) + \\
&+\Big(1-\frac{2Gm}{c^2r}\Big)(1+\bar{P}/c^2\bar{\rho})^{-1}\cdot 2y\Big)
\frac{1}{r\bar{\rho}}\frac{d\bar{P}}{dr}; 
\end{align*}
\begin{align*}
[R5]&:=
\Big(1+\frac{1}{\gamma^P}
\partial_z\Phi^P\Big)
\Big(1+\frac{r^2v^2}{c^2}-\frac{2Gm}{c^2r(1+y)}\Big)
(1+P/c^2\rho)^{-1}
\frac{(1+y)^2}{r\bar{\rho}}
\frac{\partial}{\partial r}
\bar{P}\gamma^P
(3y+z) ; \\
[R6]&:=\Big(1+\frac{r^2v^2}{c^2}-
\frac{2Gm}{c^2r(1+y)}\Big)
(1+P/c^2\rho)^{-1}(1+y)^2\frac{\bar{P}}{r\bar{\rho}}\cdot [Q0]; \\
[R7]&:=\Big(1+\frac{r^2v^2}{c^2}-
\frac{2Gm}{c^2r(1+y)}\Big)
(1+P/c^2\rho)^{-1}(1+y)^2\frac{1}{r\bar{\rho}}\frac{d\bar{P}}{dr}\cdot[Q1].
\end{align*}

Then, using (\ref{TOV2}), we have
$$[R1]=[R3L]+[R3Q]+[R4L]+[R4Q]+[R5]+[R6]+[R7].
$$

Let us define $G_1$ by
\begin{equation}
1+G_1=\Big(1+\frac{1}{\gamma^P}
\partial_z\Phi^P\Big)
\frac{1+\frac{r^2v^2}{c^2}-\frac{2Gm}{c^2r(1+y)}}{1-\frac{2Gm}{c^2r}}
\frac{1+\bar{P}/c^2\bar{\rho}}{1+P/c^2\rho}
(1+y)^2.
\end{equation}
Then
$$[R5]=(1+G_1)\Big(1-\frac{2Gm}{c^2r}\Big)
(1+\bar{P}/c^2\bar{\rho})^{-1}
\frac{1}{r\bar{\rho}}
\frac{\partial}{\partial r}\bar{P}\gamma^P(3y+z)$$
and, by definition,
\begin{align*}
-e^{-2\bar{F}}(1+\bar{P}/c^2\bar{\rho})^{-1}\mathcal{L}y&= [R3L]+[R4L]+ \\
&+\Big(1-\frac{2Gm}{c^2r}\Big)
(1+\bar{P}/c^2\bar{\rho})^{-1}\frac{1}{r\bar{\rho}}
\frac{\partial}{\partial r}\bar{P}\gamma^P(3y+z) \\
&=[R3L]+[R4L]+\frac{1}{1+G_1}[R5].
\end{align*}
This implies
\begin{align*}
[R1]&=-(1+G_1)e^{-2\bar{F}}
(1+\bar{P}/c^2\bar{\rho})^{-1}\mathcal{L}y+ \\
&-G_1([R3L]+[R4L])+[R3Q]+[R4Q]+[R6]+[R7].
\end{align*}

Now, putting
\begin{subequations}
\begin{align}
H_1&:=e^{F-2\bar{F}}(1+\bar{P}/c^2\bar{\rho})^{-1}(1+G_1), \\
H_2&:=e^FG_2, \\
G_2&:=(1+G_1)([R3L]+[R4L])-[R3]-[R4] + \nonumber \\
&-[R6]-[R7]-[R2],
\end{align}
\end{subequations}
we can write
$$e^F\times(\mbox{the right-hand side of (\ref{pert2b})})=
-H_1\mathcal{L}y-H_2.
$$

The following observation will play a crucial role in the analysis of the equation as in \cite{MakinoOJM}.

\begin{Proposition}
There is an analytic function $\hat{a}$ of $1-x, y,z, v, w, y', y''$ such that
$$(\partial_zH_1)\mathcal{L}y+\partial_zH_2=(1-x)\hat{a}$$
as $x\rightarrow 1$.
\end{Proposition}

Proof.  For the sake of abbreviations, hereafter we will denote
$$Q_1\equiv Q_0$$
if there is an analytic function $\Omega(1-x, y, z, v, w, y', y'')$ such that
$$Q_1=Q_0+(1-x)\Omega. $$

We are considering
\begin{align*}
(\partial_zH_1)\mathcal{L}y+\partial_zH_2&=
(\partial_zF\cdot H_1+
e^{F-2\bar{F}}(1+\bar{P}/c^2\rho)^{-1}\partial_zG_1)\mathcal{L}y + \\
&+\partial_zF\cdot H_2+e^F\partial_zG_2.
\end{align*}
First we note that (\ref{Fu}) and(\ref{rho}) imply
$$\partial_zF=\frac{1}{c^2}\rho\frac{du}{d\rho}(1+y+z)^{-1}$$
and that
$$\rho\frac{du}{d\rho}=(\gamma-1)u(1+[u]_1),\quad
u=\bar{u}(1+[\!(x;y,z]_1)\equiv 0.
$$
( Here $[\!(x; y, z]_1$ stands for an analytic function of $x$ in a neighborhood of
$[0,1]$ and $y,z$ of a neighborhood of $(0,0)$ of the form
$\sum_{k_1+k_2\geq 1}a_{k_1k_2}(x)y^{k_1}z^{k_2}$).
Therefore
$\partial_zF\equiv 0$ and
$$(\partial_zH_1)\mathcal{L}y+\partial_zH_2\equiv e^F[S],
$$
where
\begin{align*}
[S]&=(\partial_zG_1)e^{-2\bar{F}}(1+\bar{P}/c^2\rho)^{-1}\mathcal{L}y+
\partial_zG_2 = \\
&=-(\partial_zG_1)\Big(1-\frac{2Gm}{c^2r}\Big)
(1+\bar{P}/c^2\bar{\rho})^{-1}
\frac{1}{r\bar{\rho}}\frac{\partial}{\partial r}
\bar{P}\gamma^P(3y+z) +\\
&+(1+G_1)\frac{\partial}{\partial z}([R3L]+[R4L])-\frac{\partial}{\partial z}([R3]+[R4]+[R6]+[R7]+[R2]).
\end{align*}
But, keeping in mind that
$\bar{P}/\bar{\rho}\equiv P/\rho \equiv 0$ and that
\begin{align*}
\frac{\partial P}{\partial z}&=-\rho\frac{dP}{d\rho}(1+y+z)^{-1}\equiv 0, \\
\frac{\partial}{\partial z}\Big(\frac{P}{\rho}\Big)&=
\Big(-\frac{dP}{d\rho}+\frac{P}{\rho}\Big)(1+y+z)^{-1}\equiv 0,
\end{align*}
we see
\begin{align*}
&-(\partial_zG_1)\Big(1-\frac{2Gm}{c^2r}\Big)
(1+\bar{P}/c^2\bar{\rho})^{-1}
\frac{1}{r\bar{\rho}}\frac{\partial}{\partial r}
\bar{P}\gamma^P(3y+z) \equiv \\
&\equiv -\partial_z^2\Phi^P\Big(1+\frac{r^2v^2}{c^2}-
\frac{2Gm}{c^2r(1+y)}\Big)
(1+y)^2
\frac{1}{r\bar{\rho}}\frac{d\bar{P}}{dr}(3y+z).
\end{align*}

On the other hand it is easy to see
$$\frac{\partial}{\partial z}[R3L]\equiv
\frac{\partial}{\partial z}[R4L]\equiv
\frac{\partial}{\partial z}[R3]\equiv
\frac{\partial}{\partial z}[R4]\equiv
\frac{\partial}{\partial z}[R6]\equiv
\frac{\partial}{\partial z}[R2]\equiv 0$$
and
$$\frac{\partial}{\partial z}
[R7]\equiv
-\Big(1+\frac{r^2v^2}{c^2}-\frac{2Gm}{c^2r(1+y)}\Big)
(1+y)^2
\frac{1}{r\bar{\rho}}
\frac{d\bar{P}}{dr}(3y+z)\partial_z^2\Phi^P.
$$
Hence we have $[S]\equiv 0$ so that
$$(\partial_zH_1)\mathcal{L}y+\partial_zH_2\equiv 0.$$
This was to be shown. $\square$\\

{\bf Remark.}\hspace{5mm}  Note that $\partial[R7]/\partial z \not\equiv 0$. In fact,
we have
$$\frac{1}{r\bar{\rho}}\frac{d\bar{P}}{dr}\rightarrow -\frac{K}{r_+}\not= 0$$
and
\begin{align*}
\partial_z^2\Phi^P&=
-\frac{P}{\rho}\Big(\frac{d}{d\rho}\rho\frac{dP}{d\rho}+\frac{dP}{d\rho}\Big)(1+y+z)^{-2}
(1-\gamma^P(3y+z)-\Phi^P) \\
 &\rightarrow -(\gamma+1)(1+y+z)^{-2}(1-\gamma(3y+z)-\Phi^P(0,y,z))\not=0
\end{align*}
as $x \rightarrow 1$.\\

Now, putting
\begin{equation}
J:=e^F(1+P/c^2\rho),
\end{equation}
we rewrite the system of equations (\ref{pert1})(\ref{pert2b}) as
\begin{subequations}
\begin{align}
&\frac{\partial y}{\partial t}-Jv=0, \\
&\frac{\partial v}{\partial t}+H_1\mathcal{L}y+H_2=0.
\end{align}
\end{subequations}
Here the unknown functions are $(t,x)\mapsto y, v$. \\

\section{Framework to apply the
Nash-Moser(-Hamilton) theorem}

Having fixed a time periodic solution $Y_1$ of the linearized equation,
we put
\begin{align}
y&=\varepsilon (Y_1+Y), \\
z&=r\frac{\partial y}{\partial r}=\varepsilon (Z_1+Z) \quad\mbox{with}\quad
Z_1=r\frac{\partial Y_1}{\partial r}, \\
v&=\varepsilon(V_1+V) \quad 
\mbox{with}\quad 
V_1=\frac{1}{J^o}\frac{\partial Y_1}{\partial t}.
\end{align}
Here
\begin{equation}
J^o:=J\Big|_{y=z=0}=e^{\bar{F}}(1+\bar{P}/c^2\bar{\rho}),
\end{equation}
and $Y, Z=r\partial Y/\partial r, V$ are new unknown functions. The parameter
$\varepsilon$ will be taken sufficiently small.\\

Now the system of equations turns out to be
\begin{subequations}
\begin{align}
&\frac{\partial Y}{\partial t}-JV-(\Delta J)V_1=(J-J^o)^oV_1 
\label{eqY} \\
&\frac{\partial V}{\partial t}+H_1\mathcal{L}Y+
(\Delta H_1)(\mathcal{L}Y_1)+
\frac{1}{\varepsilon}\Delta H_2= \nonumber \\
&=-\Big(H_1-\frac{1}{J^o}\Big)^o(\mathcal{L}Y_1)
-\frac{1}{\varepsilon}H_2^o, \label{eqV}
\end{align}
\end{subequations}
where
\begin{subequations}
\begin{align}
(J-J^o)^o&:=
(J-J^o)\Big|_{Y=Z=0}
=J\Big|_{y=\varepsilon Y_1, z=\varepsilon Z_1}-J^o, \\
\Delta J&:=J-J^o-(J-J^o)^o=J-J
\Big|_{y=\varepsilon Y_1, z=\varepsilon Z_1}, \\
\Big(H_1-\frac{1}{J^o}\Big)^o&:=\Big(H_1-\frac{1}{J^o}\Big)
\Big|_{Y=Z=V=0} = \nonumber \\
&=H_1\Big|_{y=\varepsilon Y_1, z=\varepsilon Z_1, v=\varepsilon V_1}
-\frac{1}{J^o}, \\
\Delta H_1 &:= H_1-\frac{1}{J^o}-
\Big(H_1-\frac{1}{J^o}\Big)^o = \nonumber \\
&=H_1-H_1\Big|_{y=\varepsilon Y_1, z=\varepsilon Z_1, v=\varepsilon V_1}, \\
H_2^o&:=H_2\Big|_{y=\varepsilon Y_1, z=\varepsilon Z_1, v=\varepsilon V_1}, \\
\Delta H_2&= H_2-H_2^o.
\end{align}
\end{subequations}

Let us introduce the vector-valued unknown function
\begin{equation}
\vec{w}=\begin{bmatrix}
Y \\
V
\end{bmatrix}.
\end{equation}
We put
\begin{equation}
\mathfrak{P}(\vec{w})=
\begin{bmatrix}
\mbox{the left-hand side of (\ref{eqY})}\\
\mbox{the left-hand side of (\ref{eqV})}
\end{bmatrix},
\end{equation} 
and
\begin{equation}
\vec{c}=\frac{1}{\varepsilon}
\begin{bmatrix}
\mbox{the right-hand side of (\ref{eqY})} \\
\mbox{the right-hand side of (\ref{eqV})}
\end{bmatrix}
\end{equation}
The equation to be solved now is
\begin{equation}
\mathfrak{P}(\vec{w})=\varepsilon\vec{c}. 
\end{equation}

We are going to apply the Nash-Moser(-Hamilton) theorem to find
$\vec{w}=\mathfrak{P}^{-1}(\varepsilon \vec{c})$. To do it, we must analyze the Fr\'{e}chet derivative $D\mathfrak{P}$ of
the mapping $\mathfrak{P}$ at a given fixed $\vec{w}\in 
C^{\infty}([0,T]_t\times [0,1]_x)$. Introducing the new variable
\begin{equation}
\vec{h}=
\begin{bmatrix}
h \\
 k
\end{bmatrix},
\end{equation}
the Fr\'{e}chet derivative is defined by
\begin{align}
D\mathfrak{P}(\vec{w})\vec{h}&=
\lim_{s\rightarrow 0}
\frac{1}{s}(
\mathfrak{P}(\vec{w}+s\vec{h})-\mathfrak{P}(\vec{w})) \nonumber \\
&=\begin{bmatrix}
[DP1] \\
[DP2]
\end{bmatrix},
\end{align}
where
\begin{subequations}
\begin{align}
[DP1]&=\frac{\partial}{\partial t}h - J\cdot  k + \nonumber \\
&-\Big(
(\partial_yJ)v+(\partial_zJ)vr\frac{\partial}{\partial r}\Big)h, \\
[DP2]&=\frac{\partial}{\partial t} k+ H_1\cdot \mathcal{L}h + \nonumber \\
&+\Big(
(\partial_yH_1)\mathcal{L}y+\partial_yH_2+
((\partial_zH_1)\mathcal{L}y+
\partial_zH_2)r\frac{\partial}{\partial r}\Big)h + \nonumber \\
&+\Big(
(\partial_vH_1)\mathcal{L}y+\partial_vH_2+
\partial_wH_2\cdot r\frac{\partial}{\partial r}\Big) k.
\end{align}
\end{subequations}

Thanks to Proposition 11, we can claim the following

\begin{Proposition}
We have
\begin{subequations}
\begin{align}
(\partial_zJ)r\frac{\partial}{\partial r}&=[\!(x;y, Dy, D^2y, v, Dv]_0\cdot x(1-x)\frac{\partial}{\partial x}, \label{cr1} \\
((\partial_zH_1)\mathcal{L}y+
\partial_zH_2)r\frac{\partial}{\partial r}&=
[\!(x;y, Dy, D^2y,v,  Dv]_0\cdot x(1-x)\frac{\partial}{\partial x}, \label{cr2} \\
\partial_wH_2\cdot r\frac{\partial}{\partial r}&=
[\!(x; y, Dy, D^2y, v, Dv]_0\cdot x(1-x)\frac{\partial}{\partial x}.\label{cr3}
\end{align}
\end{subequations}
Here $D=\partial/\partial x$.
\end{Proposition} 

Proof.  Since
$$r\frac{\partial}{\partial r}=2x(1+[x]_1)\frac{\partial}{\partial x}$$
as
$x \rightarrow 0 \quad(r \rightarrow 0)$, the problem is concentrated to the situation as $x \rightarrow 1\quad (r\rightarrow r_+)$. Now, since 
$$\partial_zJ=(\partial_z F)J+e^F\frac{1}{c^2}\frac{\partial}{\partial z}\Big(
\frac{P}{\rho}\Big)=\frac{1}{c^2}e^F\frac{P}{\rho}(1+y+z)^{-1}, $$
it is clear that $\partial_zJ\equiv 0 \quad (\mbox{mod} (1-x) )$, that is, (\ref{cr1}).
(\ref{cr2}) is the result of Proposition 11. By definition we have
$$
\partial_wH_2=e^F\partial_wG_2=-e^F\partial_w[R2]= 
=-e^F\frac{1}{c^2}(1+y)^2\frac{P}{\bar{\rho}}v \equiv 0,
$$
that is, (\ref{cr3}). $\square$\\

In the sequel we can claim that there are analytic functions
$a_{01}, a_{00}, a_{11}, a_{10},\\
 a_{21}, a_{20}$ of
$x,y, Dy, D^2y, v, Dv$, where
$D=\partial/\partial x, y=\varepsilon (Y_1+Y), v=\varepsilon (V_1+V)$, such that the components of $D\mathfrak{P}(\vec{w})\vec{h}$ can be written as
\begin{subequations}
\begin{align}
[DP1]&=\frac{\partial}{\partial t}h-J k+
(a_{01}x(1-x)D+a_{00})h, \\
[DP2]&=\frac{\partial}{\partial t} k + H_1\mathcal{L}h +
(a_{11}x(1-x)D+a_{10})h+
(a_{21}x(1-x)D+a_{20}) k.
\end{align}
\end{subequations}

We note that $a_{01}, \cdots, a_{20} = O(\varepsilon)$ provided that
$Y, DY, V, DV=O(1)$. 

On the other hand we note, by definition, that
$$J=e^F(1+P/c^2\rho)=e^{\bar{F}}(1+\bar{P}/c^2\bar{\rho})
(1+[\!(x; y, Dy]_1)$$
and
\begin{align*}
H_1&=e^{F-2\bar{F}}(1+\bar{P}/c^2\bar{\rho})^{-1}(1+G_1) \\
&=e^{-\bar{F}}(1+\bar{P}/c^2\bar{\rho})^{-1}
(1+[\!(x;y,Dy]_1+
v^2[\!(x;y,Dy]_0).
\end{align*}
Hence $$J=e^{\bar{F}}(1+\bar{P}/c^2\bar{\rho})(1+O(\varepsilon)),$$
$$H_1=e^{-\bar{F}}(1+\bar{P}/c^2\bar{\rho})^{-1}(1+O(\varepsilon))$$
provided that $Y, DY=O(1)$.\\

{\bf Remark.} \hspace{5mm} We see $J\rightarrow 1, H_1
 \rightarrow 1+[y,z]_1$ as $1/c^2
\rightarrow 0$, while $\bar{P}/\bar{\rho}, \bar{u}$ are supposed
to be bounded. (The equilibrium depends upon the central density 
$\rho_c$ and the speed of light $c$.)
 But we do not discuss the details of the non-relativistic limit in this article. 

\section{Main conclusion}

Now we are ready to propose the main conclusion of this article.

\begin{Theorem}
Given $T>0$, there is a positive number $\varepsilon_0(T)$ such that, for
$|\varepsilon|\leq\varepsilon_0(T)$, there is a solution 
$\vec{w}\in C^{\infty}([0,T]\times[0,1])$ of  (7.10)  such that
$$\sup_{j+k\leq n}\Big\|\Big(\frac{\partial}{\partial t}\Big)^j
\Big(\frac{\partial}{\partial x}\Big)^k\vec{w}
\Big\|_{L^{\infty}([0,T]\times[0,1])}\leq C_n|\varepsilon|,
$$
and hence a solution $(y,v)\in C^{\infty}([0,T]\times[0,r_+])$ of (\ref{pert1})(\ref{pert2b}) of the form
$$y=\varepsilon Y_1+O(\varepsilon^2).$$
\end{Theorem}

Note that for this solution the component $R$ of the metric (\ref{metric}) behaves like
$$R=r(1+\varepsilon Y_1+O(\varepsilon^2)),$$
and the density distribution enjoys
$$\rho=\begin{cases}
C(t)(r_+-r)^{1/(\gamma-1)}(1+O(r_+-r)) & (0\leq r<r_+) \\
0 & (r+\leq r)
\end{cases}
$$
Here $C(t)$ is a smooth positive function of $t$.

In other words, the value $R_+(t)$ of the Eulerian coordinate $R$ at the surface of the star
$r=r_+$ is approximately oscillating as
$$R_+(t)=r_+(1+\varepsilon
\sin(\sqrt{\lambda}t+\Theta_0)\psi(1)+
O(\varepsilon^2)).$$\\

A proof can be given by an application of the Nash-Moser(-Hamilton)
theorem (\cite[p.171, III.1.1.1.]{Hamilton}) as in \cite{MakinoOJM}, \cite{MakinoFE}. The discussion is quite parallel. Therefore, omitting the repetitions of the details, we will explain only the points for which some modifications are necessary.\\

First the mapping $\mathfrak{P}$ is considered on the tame spaces
$\vec{\mathfrak{E}}$ and $\vec{\mathfrak{E}}_0$. Here
$\vec{\mathfrak{E}}=\mathfrak{E}\times\mathfrak{E}$ with 
$\mathfrak{E}=C^{\infty}([0,T]\times[0,1])$ and
$\vec{\mathfrak{E}}_0=\mathfrak{E}_0\times\mathfrak{E}_0$ with
$\mathfrak{E}_0=\{ \phi\in
\mathfrak{E} \quad | \quad \phi=0 \quad \mbox{at}\quad t=0\}$. 
Since $\mathfrak{E}$ admits the gradings of norms as in \cite{MakinoOJM},
$\vec{\mathfrak{E}}$ is a tame space as the direct Cartesian
product. The domain of $\mathfrak{P}$ is $\vec{\mathfrak{U}}$, the set of
all functions $\vec{w}=(Y,V)^T\in\vec{\mathfrak{E}}_0$ such that
$$|Y|+|DY|+|V|+|DV|<1.$$
We consider $\varepsilon$ such that $|\varepsilon|\leq \varepsilon_1$, $\varepsilon_1$ being a fixed sufficiently small positive number. The mapping $\mathfrak{P}$ is a tame mapping from $\vec{\mathfrak{U}}$ into
$\vec{\mathfrak{E}}$.\\

Introducing the operator
\begin{equation}
\Lambda=x(1-x)\frac{d^2}{dx^2}+\Big(\frac{5}{2}(1-x)-
\frac{N}{2}x\Big)\frac{d}{dx},
\end{equation}
just as \cite[(20)]{MakinoOJM}, we rewrite
the second component of
$D\mathfrak{P}(\vec{w})\vec{h}$ as
\begin{align}
[DP2]&=\frac{\partial}{\partial t} k-H_1\Lambda h + \nonumber \\
&+b_1\check{D}h+ b_0h+
a_{21}\check{D} k+
a_{20} k,
\end{align}
where
\begin{subequations}
\begin{align}
\check{D}&=x(1-x)\frac{\partial}{\partial x}, \\
b_1&=\frac{H_1L_1}{x(1-x)}+a_{11}, \\
b_0&=H_1L_0+a_{10}.
\end{align}
\end{subequations}
Then $b_1, b_0, a_{21}, a_{20} $ are analytic functions of
$x,y,Dy, D^2y, v, Dv$.\\

Let us introduce the Hilbert spaces $\mathfrak{X}=\mathfrak{X}^0,
\mathfrak{X}^1, \mathfrak{X}^2$, just in the same manner as \cite{MakinoOJM}, by
\begin{align*}
\mathfrak{X}&=L^2((0,1); x^{\frac{3}{2}}(1-x)^{\frac{N}{2}-1}dx), \\
\mathfrak{X}^1&=\{ \phi\in \mathfrak{X}\quad|\quad \dot{D}\phi:=
\sqrt{x(1-x)}\frac{d\phi}{dx}\in\mathfrak{X}\},\\
\mathfrak{X}^2&=\{\phi\in\mathfrak{X}^1\quad|\quad -\Lambda \phi \in \mathfrak{X}\}.
\end{align*}\\

We write the equation
$$D\mathfrak{P}(\vec{w})\vec{h}=\vec{g}, $$
where $\vec{g}=(g_1,g_2)^T$ is a given function in $\vec{\mathfrak{E}}$, as
\begin{equation}
\frac{\partial}{\partial t}
\begin{bmatrix}
h \\
 k
\end{bmatrix}
+
\begin{bmatrix}
\mathfrak{a}_1 & -J \\
\mathcal{A} & \mathfrak{a}_2
\end{bmatrix}
\begin{bmatrix}
h \\
 k
\end{bmatrix}
=
\begin{bmatrix}
g_1 \\
g_2
\end{bmatrix},
\label{EqLvec}
\end{equation}
where
\begin{subequations}
\begin{align}
\mathfrak{a}_1&=a_{01}\check{D}+a_{00}, \\
\mathfrak{a}_2&=a_{21}\check{D}+a_{20,} \\
\mathcal{A}&=-H_1\Lambda + b_1\check{D}+b_0.
\end{align}
\end{subequations}

Then the standard calculation leads us to the equality
\begin{align*}
&\frac{1}{2}\frac{d}{dt}\Big[\| k\|^2+\Big(
\frac{H_1}{J}\dot{D}h\Big|\dot{D}h\Big)\Big] + \\
&+(\beta_1\dot{D}h|\dot{D}h)+(\beta_2\dot{D}h|h) + (\beta_3\dot{D}h| k) 
+(\beta_4h| k)+ (\beta_5 k| k) = \\
&=\Big(\frac{H_1}{J}\dot{D}h\Big|\dot{D}g_1\Big)+( k|g_2),
\end{align*}
where
\begin{align*}
\beta_1&=-\frac{1}{4}(3+(N+3)x+2\check{D})\frac{H_1a_{01}}{J}
-\frac{1}{2}\frac{\partial}{\partial t}\frac{H_1}{J}+
\frac{H_1}{J}(\check{D}a_{01}+a_{00}), \\
\beta_2&=\frac{H_1}{J}\dot{D}a_{00}, \\
\beta_3&=-\frac{H_1}{J}\dot{D}J+\dot{D}H_1+
\sqrt{x(1-x)}(b_1+a_{21}),\\
\beta_4&=b_0, \\
\beta_5&=a_{20}.
\end{align*}
Here
$$(\phi|\psi)=(\phi|\psi)_{\mathfrak{X}}=\int_0^1
\phi\bar{\psi}x^{\frac{3}{2}}(1-x)^{\frac{N}{2}-1}dx
\quad \mbox{and}\quad
\|\phi\|=\|\phi\|_{\mathfrak{X}}=\sqrt{(\phi|\phi)_{\mathfrak{X}}},$$
and we have used the formula
$$(\alpha\dot{D}h|\dot{D}\check{D}h)=(\alpha^*\dot{D}h|\dot{D}h)\quad\mbox{with}\quad
\alpha^*=-\frac{1}{4}(3+(N+3)x+2\check{D})\alpha,
$$
which holds for $h \in \mathfrak{X}^2$ and $\alpha\in C^{\infty}([0,1])$,
together with \cite[Proposition 8]{MakinoOJM}.

Since $\vec{w}$ is confined to $\vec{\mathfrak{U}}$ and $|\varepsilon|$
is restricted $\leq\varepsilon_0$, we can assume
$$\frac{1}{M_0}\leq J\leq M_0, \qquad \frac{1}{M_0}\leq H_1\leq M_0$$
with a constant $M_0$ independent of $\vec{w}$.\\

Now 
the energy
$$\mathcal{E}:=\| k\|^2+\Big(
\frac{H_1}{J}\dot{D}h|\dot{D}h\Big)$$
enjoys the inequality
$$\frac{1}{2}\frac{d\mathcal{E}}{dt}\leq M\Big(\|\vec{h}
\|_{\mathfrak{H}}^2+
\|\vec{h}\|_{\mathfrak{H}}\|\vec{g}\|_{\mathfrak{H}}\Big), $$
where
$\mathfrak{H}=\mathfrak{X}^1\times\mathfrak{X}$ and 
$$
\|(\phi,\psi)^T\|_{\mathfrak{H}}^2=
\|\phi\|_{\mathfrak{X}^1}^2+\|\psi\|_{\mathfrak{X}}^2
=\|\phi\|^2+\|\dot{D}\phi\|^2+\|\psi\|^2,$$
and
$$M=\sum_{1\leq j\leq 5}\|\beta_j\|_{L^{\infty}}+ (M_0)^2+1.
$$
Since
$$\frac{1}{(M_0)^2}(\| k\|^2+
\|\dot{D}h\|^2)
\leq\mathcal{E}
\leq
(M_0)^2(\| k\|^2+
\|\dot{D}h\|^2),
$$
using the same Gronwall's argument as \cite[Proposition 9]{MakinoOJM},
\cite[Lemma 3]{MakinoFE}, we see that the initial value problem for the equation
(\ref{EqLvec}) with the initial condition
$$h= k=0 \qquad \mbox{at}\qquad t=0$$
admits a unique solution
$\vec{h}=(h,  k)^T$ in
$C([0,T], \mathfrak{X}^2\times\mathfrak{X}^1)$ for given
$\vec{g}\in C([0,T], \mathfrak{X}^1\times\mathfrak{X})$, which enjoys the energy estimate
$$\|\vec{h}\|_{\mathfrak{H}}\leq
C\int_0^t
\|\vec{g}(t')\|_{\mathfrak{H}}dt'.$$
Therefore $D\mathfrak{P}(\vec{w})$ admits an inverse, and its tame estimates can be shown in the same manner as \cite{MakinoOJM}. 
An outline of this procedure can be found in Appendix. This completes the proof of the main conclusion.

\section{Cauchy problems}

As a supplement let us consider the Cauchy problem
associated with equations (6.10a)(6.10b), that is, (CP):
\begin{align}
&\frac{\partial y}{\partial t}-Jv=0, \quad
\frac{\partial v}{\partial t}+H_1\mathcal{L}y+H_2=0,\quad(t\geq 0, 0\leq x\leq 1) \\
&y|_{t=0}=\psi_0(x),\quad v|_{t=0}=\psi_1(x).
\end{align}
Here $\psi_0$ and $\psi_1$ are functions given in $C^{\infty}([0,1])$.\\

Let us recall that
$$J=e^F(1+P/c^2\rho)=J(x,y,z)$$
is an analytic function of $x$ (in a neighborhood of $[0,1]$),
$y$ (small) and $\displaystyle z=r\frac{\partial y}{\partial r}=
x[\!(x )\!]\frac{\partial y}{\partial x}$ (small),
where $[\!(x)\!]$ stands for an analytic function of $x$ in a neighborhood of $[0,1]$; that $H_1$ and 
$H_2$ are analytic functions of $x,y,z,v$ and $\displaystyle w=r
\frac{\partial v}{\partial r}$ (quadratic in $v/c, w/c$); and that
the linear operator $\mathcal{L}$ has the form
\begin{align*}\mathcal{L}y=&
-x(1-x)\frac{d^2y}{dx^2}-
\Big(\frac{5}{2}(1-x)-\frac{N}{2}x\Big)\frac{dy}{dx}
+\\
&+l_1(x)x(1-x)\frac{dy}{dx}+L_0(x)y,
\end{align*}
where $l_1$ and $L_0$ are analytic functions of $x$ in a neighborhood of 
$[0,1]$.\\

We claim the following
\begin{Theorem}
For any given $T>0$ there exists a sufficiently small positive number $\delta$ such that
if $\psi_0, \psi_1\in C^{\infty}([0,1])$ satisfy
$$\max_{k\leq\mathfrak{K}}
\Big\{
\Big\|\Big(
\frac{d}{dx}
\Big)^k\psi_0\Big\|_{L^{\infty}},
\Big\|\Big(
\frac{d}{dx}
\Big)^k\psi_1\Big\|_{L^{\infty}}
\Big\}\leq\delta $$
then there exists a unique solution $(y,v)$ of (CP)
in $C^{\infty}([0,T]\times[0,1])$. Here $\mathfrak{K}$ is a 
sufficiently large number depending only upon $\gamma$.
\end{Theorem}

Proof can be done in a almost same manner as that of Theorem 2 of
\cite{MakinoOJM}. Let us take the functions
\begin{equation}
y_1^*=\psi_0(x)+tJ^o(x)\psi_1(x),\quad
v_1^*=\psi_1(x),
\end{equation}
which satisfy the initial conditions, where
$$J^o(x)=J(x, 0, 0)$$
as (7.4). Then we seek a solution $(y,v)$ of the form
\begin{equation}
y=y_1^*+Y, \qquad
v=v_1^*+V.
\end{equation}
The initial condition for 
$\vec{w}:=(Y,V)^T$ is
\begin{equation}
\vec{w}|_{t=0}=(0,0)^T
\end{equation}
and the equations to be satisfied by $\vec{w}=(Y,V)^T$ are
\begin{subequations}
\begin{align}
&\frac{\partial Y}{\partial t}-JV-
(\Delta J)v_1^*=c_1, \\
&\frac{\partial V}{\partial t}+
H_1\mathcal{L}Y+
(\Delta H_1)\mathcal{L}y_1^*+
\Delta H_2=c_2,
\end{align}
\end{subequations}
where
\begin{subequations}
\begin{align}
J&=J(x,y_1^*+Y, z_1^*+Z),\quad \mbox{with}\quad Z:=r
\frac{\partial Y}{\partial r}, \\
\Delta J&=J(x,y_1^*+Y, z_1^*+Z)-J(x,y_1^*, z_1^*), \\
c_1&=J(x,y_1^*,z_1^*)-J(x,0,0), \\
H_1&=H_1(x,y_1^*+Y, z_1^*+Z, v_1^*+V, w_1^*+W),\quad
\mbox{with}\quad W=r\frac{\partial V}{\partial r}, \\
\Delta H_1&=H_1(x,y_1^*+Y, z_1^*+Z, v_1^*+V, w_1^*+W)
-H_1(x,y_1^*, z_1^*, v_1^*, w_1^*)\\
\Delta H_2&==H_2(x,y_1^*+Y, z_1^*+Z, v_1^*+V, w_1^*+W)
-H_2(x,y_1^*, z_1^*, v_1^*, w_1^*)\\
c_2&=-H_1(x,y_1^*,z_1^*,v_1^*,w_1^*)\mathcal{L}y_1^*
-H_2(x,y_1^*,z_1^*,v_1^*, w_1^*).
\end{align}
\end{subequations}
The problem can be written as
\begin{equation}
\mathfrak{P}(\vec{w})=\vec{c},
\end{equation}
where
\begin{align}
\mathfrak{P}(\vec{w})&=
\begin{bmatrix}
\mbox{the left-hand side of (9.6a)}\\
\mbox{the left-hand side of (9.6b)}
\end{bmatrix}\\
\vec{c}&=\begin{bmatrix}
c_1 \\
c_2
\end{bmatrix}.
\end{align}
Then the Nash-Moser(-Hamilton) theorem
can be applied in the same manner as the proof of
Theorem 1, since the Fr\'{e}chet derivative of
$\mathfrak{P}$ has the same form as (7.12)(7.13). This completes the proof.\\

{\bf Remark.} The initial data read
\begin{align*}
R|_{t=0}&=r(1+\psi_0(x(r))), \\
\frac{\partial R}{\partial t}\Big|_{t=0}&=
\frac{1}{c}\sqrt{1-\frac{2Gm_+}{c^2r_+}}
\exp\Big[-\frac{1}{c^2}u(\rho^0)\Big]r\psi_1(x(r)),
\end{align*}
where
$$\rho^0=\bar{\rho}(r)
(1+\psi_0)^{-2}
\Big(1+\psi_0+
r\frac{d\psi_0}{dr}\Big)^{-1}.
$$

{\bf Supplementary Remark 4}.
Let us consider moving solutions constructed in Section 8 or 9,
which are defined on $0\leq t\leq T, 0\leq r\leq r_+$. 
We should discuss how to extend the metric onto the exterior vacuum region $r>r_+$.
We owe the idea to \cite{MisnerS}

If a spherically symmetric extension to the vacuum region is possible, the Birkhoff's theorem reads that it should be the Schwartzschild's metric
$$ds^2=\Big(1-\frac{2Gm_+}{c^2R^{\sharp}}\Big)c^2(dt^{\sharp})^2
-\Big(1-\frac{2Gm_+}{c^2R^{\sharp}}\Big)^{-1}(dR^{\sharp})^2-
(R^{\sharp})^2(d\theta^2+
\sin^2\theta d\phi^2).$$
Here $t^{\sharp}=t^{\sharp}(t,r), R^{\sharp}=R^{\sharp}(t,r)$ are smooth functions of
$0\leq t\leq T, r_+\leq r<\infty$. We have:

{\bf  There are $t^{\sharp}(t,r), R^{\sharp}(t,r)$ such that the components of the metric are of class $C^1([0,T]\times[0, +\infty))$.}

Let us verify this. We are considering the patched metric
$$ds^2=g_{00}c^2dt^2+
2g_{01}cdtdr+
g_{11}dr^2+
g_{22}(d\theta^2+\sin^2\theta d\phi^2),$$
where
\begin{align*}
g_{00}&=\begin{cases}
\kappa e^{-2u/c^2}  \quad (0\leq r\leq r_+) \\
\displaystyle K^{\sharp}\Big(\frac{\partial t^{\sharp}}{\partial t}\Big)^2-
\frac{1}{c^2}(K^{\sharp})^{-1}
\Big(\frac{\partial R^{\sharp}}{\partial t}\Big)^2  \quad (r_+<r),
\end{cases}\\
g_{01}&=\begin{cases} 
0  \quad (0\leq r\leq r_+) \\
\displaystyle cK^{\sharp}\frac{\partial t^{\sharp}}{\partial t}
\frac{\partial t^{\sharp}}{\partial r}-\frac{1}{c}
(K^{\sharp})^{-1}\frac{\partial R^{\sharp}}{\partial t}\frac{\partial R^{\sharp}}{\partial r}
\quad (r_+<r),
\end{cases}
\\
g_{11}&=\begin{cases}
\displaystyle -\Big(1+\frac{V^2}{c^2}-
\frac{2Gm}{c^2R}\Big)^{-1}\Big(\frac{\partial R}{\partial r}\Big)^2 
\quad (0\leq r\leq r_+) \\
\\
\displaystyle c^2K^{\sharp}\Big(\frac{\partial t^{\sharp}}{\partial r}\Big)^2
-(K^{\sharp})^{-1}
\Big(\frac{\partial R^{\sharp}}{\partial r}\Big)^2  \quad (
r_+
<r)
\end{cases}\\
g_{22}&=\begin{cases}
-R^2  \quad (0\leq t\leq r_+) \\
-(R^{\sharp})^2 \quad (r_+<r).
\end{cases}
\end{align*}
Here
$$K^{\sharp}=1-\frac{2Gm_+}{c^2R^{\sharp}}.$$

Let us assume $R=R^{\sharp}, \partial_r R=\partial_r R^{\sharp}$ at $r=r_+$
in order that $g_{22}$ be of class $C^1$.

First in order that $g_{00}$ be continuous across $r=r_+$ we require
\begin{equation}
\frac{\partial t^{\sharp}}{\partial t}=
\sqrt{\kappa}(K^{\sharp})^{-1}\Big(1+\frac{V^2}{c^2}-
\frac{2Gm_+}{c^2R}\Big)^{1/2}, \label{Eq8}
\end{equation}
on $r=r_+$, where $V=V(t, r_+-0) (=\displaystyle \frac{1}{\sqrt{\kappa}}\frac{\partial R}{\partial t})$.
In order that $g_{01}$ be continuous,
we require
\begin{equation}
\frac{\partial t^{\sharp}}{\partial r}=
\frac{1}{c^2}(K^{\sharp})^{-1}\Big(1+\frac{V^2}{c^2}-
\frac{2Gm_+}{c^2R}\Big)^{-1/2}
V
\frac{\partial R}{\partial r} \label{Eq9}
\end{equation}
on $r=r_+$. It can be shown that (\ref{Eq9}) is sufficient in order
that $g_{11}$ be continuous across $r=r_+$.
Summing up, $g_{\mu\nu}$ are continuous if (\ref{Eq8}) and 
(\ref{Eq9}) hold.
Note that, since 
$$K^{\sharp}\doteqdot \kappa,\quad
1+\frac{V^2}{c^2}-\frac{Gm_+}{c^2R}\doteqdot \kappa,\quad \frac{\partial R}{\partial r}\doteqdot 1,$$
the right-hand side of (\ref{Eq9})$\doteqdot V/c^2$ so that
$\partial t^{\sharp}/\partial r \not= 0$ and $t^{\sharp}$ should actually
depend upon $r$ if $V\not= 0$, that is, if the solution is  actually moving.

By a tedious calculation we can show that the differentiation of (\ref{Eq8})
with respect to $t$ gives the continuity of $\partial_rg_{00}$. 
On the other hand the continuity of $\partial_rg_{01}$ reads a condition of the form
\begin{equation}
K^{\sharp}\frac{\partial t^{\sharp}}{\partial t}\frac{\partial^2t^{\sharp}}{\partial r^2}-
\frac{1}{c^2}(K^{\sharp})^{-1}
\frac{\partial R}{\partial t}\frac{\partial^2R^{\sharp}}{\partial r^2}=b_1
\label{Eq10}
\end{equation}
on $r=r_+$, where $b_1$ is a function of the values of
$\partial_t t^{\sharp}, \partial_r t^{\sharp},
\partial_t\partial_rt^{\sharp}, R, \partial_rR, \partial_t\partial_rR$ on
$r=r_+$.
The continuity of $\partial_rg_{11}$ reads a condition of the form
\begin{equation}
c^2K^{\sharp}\frac{\partial t^{\sharp}}{\partial r}\frac{\partial^2t^{\sharp}}{\partial r^2}
-(K^{\sharp})^{-1}\frac{\partial R}{\partial r}
\frac{\partial^2R^{\sharp}}{\partial r^2}=b_2,\label{Eq11}
\end{equation}
on $r=r_+$, where $b_2$ is a function of the same kind as $b_1$.
If we consider (\ref{Eq10})(\ref{Eq11}) as a system of  simultaneous linear equations for the unknown
${\partial^2t^{\sharp}}/{\partial r^2}, 
{\partial^2R^{\sharp}}/{\partial r^2}$, the determinant of
the coefficient matrix is
$$-\sqrt{\kappa}\Big(1+\frac{V^2}{c^2}-\frac{2Gm_+}{c^2R}\Big)^{-1/2}
\frac{\partial R}{\partial r},$$
which is near to $-1$, since $$1+\frac{V^2}{c^2}-
\frac{2Gm_+}{c^2R} \doteqdot \kappa, \qquad \frac{\partial R}{\partial r} \doteqdot 1.$$
Since $b_1, b_2$ are known by (\ref{Eq8})(\ref{Eq9}), the values 
$\partial^2t^{\sharp}/\partial r^2, \partial^2R^{\sharp}/\partial r^2$ along $r=r_++0$ are uniquely determined.
Then all $g_{\mu\nu}$ are of class $C^1$. 

However we note that this
$\partial^2R^{\sharp}/\partial r^2$ 
generally does not coincide with $\partial^2R/\partial r^2$ on $r=r_+$,  which coincidence is necessary to that $g_{22}$ be twice continuously differentiable. 
In fact by a tedious calculation
we get 
$$\frac{\partial^2R^{\sharp}}{\partial r^2}\Big|_{r=r_++0}=
\mathcal{A}\Big(\frac{\partial R}{\partial r}\Big)^2+
\frac{\partial^2R}{\partial r^2}\Big|_{r=r_+-0},$$
where
$$\mathcal{A}=-\frac{V^2}{c^2}
\Big(
\frac{Gm_+}{c^2R^2}+
\frac{1}{\sqrt{\kappa}}
\frac{1}{c^2}\frac{\partial V}{\partial t}\Big)
\Big(1+\frac{V^2}{c^2}-
\frac{2Gm_+}{c^2R}\Big)^{-2}$$
evaluated at $r=r_+-0$. Since
$$\frac{Gm_+}{c^2R^2}\Big(1+\frac{V^2}{c^2}-
\frac{2Gm_+}{c^2R}\Big)^{-2}\doteqdot 
\frac{Gm_+}{c^2r_+^2\kappa^2}\not= 0,$$
we see that $\partial^2R^{\sharp}/\partial r^2\equiv
\partial^2R/\partial r^2$ if and only if
$V \equiv 0$
at $r=r_+-0$, which is the case
if the solution 
under consideration is an equilibrium.

Anyway we have determined the functions
\begin{align*}
f_0(t)&:=R(t,r_+),\qquad f_1(t):=\partial_rR(t,r_+), \\
f_2(t)&:=\frac{\partial^2R^{\sharp}}{\partial r^2}\quad\mbox{at}\quad
r=r_++0, \\
H(t)&:=\frac{\partial t^{\sharp}}{\partial t}\quad\mbox{at}\quad
r=R_++0,\qquad h_0(t):=
\int_0^tH(t')dt', \\
h_1(t)&:=\frac{\partial t^{\sharp}}{\partial r}\quad\mbox{at}
\quad r=r_++0,\qquad
h_2(t):=\frac{\partial^2t^{\sharp}}{\partial r^2}
\quad\mbox{at}\quad r=r_++0
\end{align*}
for $0\leq t\leq T$. Using these functions we define
$t^{\sharp}(t,r), R^{\sharp}(t,r)$ for $0\leq t\leq T,
r_+\leq r <+\infty$ as follow:
\begin{align*}
R^{\sharp}(t,r)&=f_0(t)+
f_1(t)(r-r_+)+
\frac{1}{2}f_2(t)(r-r_+)^2\chi(r-r_+), \\
t^{\sharp}(t,r)&=h_0(t)+
\Big(h_1(t)(r-r_+)+
\frac{1}{2}h_2(t)(r-r_+)^2\Big)\chi(\delta(r-r_+)).
\end{align*}
Here $\chi$ is a smooth cut off function $\in
C^{\infty}[0,+\infty)$ such that
$0\leq \chi(s)\leq 1, \chi(s)=1 $ for
$0\leq s\leq 1$ and $\chi(s)=0$ for
$2\leq s <+\infty$ and $\delta$ is a sufficiently small positive
number.
Since $f_0(t)\doteqdot r_+, f_1(t)\doteqdot 1, f_2(t)\doteqdot 0,
H(t)\doteqdot 1$, we see that
$\partial R^{\sharp}/\partial r \doteqdot 1$ and
$\partial t^{\sharp}/\partial t\doteqdot 1$ uniformly.
Then the coefficients of the metric $g_{00}, g_{01}, g_{11} $ and
$g_{22}$ are of class $C^1([0,T]\times[0, +\infty))$ and
their second order derivatives may have discontinuity of at most the first kind along the segment
$r=r_+$, and satisfy the Einstein equations in the usual sense on 
$r\not=r_+$. So, we can say that this metric is a weak solution of
the Einstein equations on $[0,T]\times \mathbb{R}^3$
in the following sense:
The Einstein equations can be written as
$$R_{\mu\nu}=\frac{8\pi G}{c^4}\Big(T_{\mu\nu}-\frac{1}{2}g_{\mu\nu}T\Big)$$
and 
\begin{align*}
T&=T^{\alpha\beta}T_{\alpha\beta}, \\
R_{\mu\nu}&=\frac{1}{2}g^{\alpha\beta}(
-\partial_{\alpha}\partial_{\beta}g_{\mu\nu}
-\partial_{\mu}\partial_{\nu}g_{\alpha\beta}+
\partial_{\beta}\partial_{\nu}g_{\mu\alpha}
+\partial_{\mu}\partial_{\alpha}g_{\beta\nu})+
F_{\mu\nu}, \\
F_{\mu\nu}&=\frac{1}{2}\partial_{\alpha}g^{\alpha\beta}(
\partial_{\nu}g_{\beta\mu}+\partial_{\mu}g_{\beta\nu}-
\partial_{\beta}g_{\mu\nu}) + \\
&-\frac{1}{2}\partial_{\nu}g^{\alpha\beta}
(\partial_{\alpha}g_{\beta\mu}+
\partial_{\mu}g_{\beta\alpha}-
\partial_{\beta}g_{\mu\alpha});
\end{align*}
Therefore, $(T_{\mu\nu})_{\mu\nu}, T \in L^2_{{loc}}$ given,
$ds^2=g_{\mu\nu}dx^{\mu}dx^{\nu}$ is said to be a weak
solution if $g_{\mu\nu}, g^{\alpha\beta} \in H^1_{{loc}}$ and for
any
test function $(\phi^{\mu\nu})_{\mu\nu}$ there holds
\begin{align*}
&\frac{1}{2}\int(\partial_{\beta}g_{\mu\nu})\partial_{\alpha}
(g^{\alpha\beta} \phi^{\mu\nu})+
(\partial_{\nu}g_{\alpha\beta})
\partial_{\mu}(g^{\alpha\beta}\phi^{\mu\nu})+ \\
&-(\partial_{\nu}g_{\mu\alpha})\partial_{\beta}(g^{\alpha\beta}\phi^{\mu\nu})
-(\partial_{\alpha}g_{\beta\nu})\partial_{\mu}(g^{\alpha\beta}\phi^{\mu\nu})
+\int F_{\mu\nu}\phi^{\mu\nu} = \\
&=\frac{8\pi G}{c^4}\int\Big(T_{\mu\nu}-\frac{1}{2}g_{\mu\nu}T\Big)\phi^{\mu\nu}.
\end{align*}
 
\medskip

{\bf\large Acknowledgment} A part of this work was done during the stay of the author at Universit\'{e} de Strasbourg at the opportunity of `94e Rencontre entre math\'{e}maticiens et physiciens th\'{e}oriciens: Riemann, Einstein et la geom\'{e}trie'. The author expresses his thanks to the organizers Professor Athanase Papadopoulos (Strasbourg) and Professor Sumio Yamada (Tokyo) for their hospitality. 

The author expresses his sincere thanks to the anonymous referee
for his/her careful reading of the original manuscript and many 
kind suggestions to ameliorate the expressions.\\

{\bf\Large Appendix}\\

Let us give an outline of the tame estimate of the mapping
$(\vec{w}, \vec{g})\mapsto \vec{h}$ when
$D\mathfrak{P}(\vec{w})\vec{h}=\vec{g}$.
The equation (\ref{EqLvec}) is split as \cite{MakinoOJM} using
a cut off function $\omega\in C^{\infty}$ such that
$\omega(x)=1$ for $x\leq 1/3$, $0<\omega(x)<1$ for
$1/3<x<2/3$ and $\omega(x)=0$ for $2/3\leq x$. Put
$$\vec{h}^{[0]}(x)=\omega(x)\vec{h}(x),\qquad
\vec{h}^{[1]}(x)=(1-\omega(x))\vec{h}(x). $$
The equations turn out to be
\begin{align*}
\frac{\partial}{\partial t}
\begin{bmatrix}
h^{[\mu]} \\
k^{[\mu]}
\end{bmatrix}
&+
\begin{bmatrix}
\mathfrak{a}_1^{[\mu]} & -J \\
\mathcal{A}^{[\mu]} & \mathfrak{a}_2^{[\mu]}
\end{bmatrix}
\begin{bmatrix}
h^{[\mu]} \\
k^{[\mu]}
\end{bmatrix} = \\ 
&=\begin{bmatrix}
g_1^{[\mu]} \\
g_2^{[\mu]}
\end{bmatrix}
+(-1)^{\mu}
\begin{bmatrix}
c_{11} & 0 \\
\mathfrak{c}_{21} & c_{22}
\end{bmatrix}
\begin{bmatrix}
h^{[1-\mu]} \\
k^{[1-\mu]}
\end{bmatrix},
\end{align*}
where $\mu=0,1$ and
\begin{align*}
\mathfrak{a}_1^{[\mu]}&=a_{01}\check{D}+a_{00}-(-1)^{\mu}a_{01}
\check{D}\omega, \\
\mathfrak{a}_2^{[\mu]}&=a_{21}\check{D}+a_{20}-(-1)^{\mu}
a_{21}\check{D}\omega, \\
\mathcal{A}^{[\mu]}&=-H_1\Lambda
+(b_1+(-1)^{\mu}2H_1(D\omega))\check{D}+
b_0+(-1)^{\mu}(H_1\Lambda-b_1\check{D})\omega, \\
c_{11}&=a_{01}\check{D}\omega, \\
\mathfrak{c}_{21}&=
-2H_1(D\omega)\check{D}+b_1(\check{D}\omega)-H_1(\Lambda\omega), \\
c_{22}&=a_{21}\check{D}\omega.
\end{align*}

Therefore the problem is reduced to the tame estimate of an equation of the form
$$\frac{\partial\vec{h}}{\partial t}+\mathfrak{A}\vec{h}=\vec{g},$$
$$\mathfrak{A}=
\begin{bmatrix}
\mathfrak{a}_1 & J \\
\mathcal{A} & \mathfrak{a}_2
\end{bmatrix}
=\begin{bmatrix}
a_{01}\check{D}+a_{00} & J \\
-b_2\triangle +b_1\check{D}+b_0 & a_{21}\check{D}+a_{20}
\end{bmatrix}
$$
under the boundary condition $h|_{x=1}=0$, where
$$ \triangle =x\frac{d^2}{dx^2}+\frac{N}{2}\frac{d}{dx},
\qquad \check{D}=x\frac{d}{dx} $$
with $N$ standing for either $2\gamma/(\gamma-1)$ or $5$.

As in \cite{MakinoOJM}, we use the notations
\begin{align*}
\vec{a}&=(a_i)_{i=0}^7=(b_0, b_1, b_2, a_{01}, a_{00}, a_{21}, a_{20}, J), \\
|\vec{a}|_n^{\langle T\rangle}&=\sup_{0\leq t\leq T}|\vec{a}|_n, \\
|\vec{a}|_n&=
\max_{j+k\leq n, 0\leq i\leq 7}
\|\partial_t^j\dot{D}^ka_i\|_{L^{\infty}} \\
\|\vec{h}\|_n^{\langle T\rangle}&=\Big(
\sum_{j+k\leq n}
\int_0^T
\|\partial_t^j\vec{h}\|_k^2dt \Big)^{1/2} \\
\|\vec{h}\|_k&=\Big(
\sum_{0\leq \ell \leq k}\langle h\rangle _{\ell+1}^2+
\langle k\rangle_{\ell}^2\Big)^{1/2}.
\end{align*}
Here $\langle \phi\rangle_{\ell}$ means the same as \cite{MakinoFE}.\\

Then the elliptic a priori estimate \cite[Proposition 8]{MakinoFE} should read
$$\|\vec{h}\|_{n+1}\leq C
(\|\mathfrak{A}\vec{h}\|_n+
(1+|\vec{a}|_{n+4})\|\vec{h}\|_0).
$$
This can be verified if we keep in mind that
\begin{align*}
\|\mathfrak{a}_1h\|_1&\leq C(|\varepsilon|\|h\|_2+\|h\|_1), \\
\|\mathfrak{a}_2k\|_0&\leq C(|\varepsilon|\|k\|_1+\|k\|_0),
\end{align*}
which come from
\begin{align*}
a_{01}&=\frac{1}{c^2}e^F\frac{P}{\rho}(1+y+z)^{-1}\varepsilon(V_1+V)
\frac{r}{x(1-x)}\frac{dx}{dr}, \\
a_{21}&=
-\frac{1}{c^2}e^F\frac{P}{\bar{\rho}}(1+y)^2\varepsilon(V_1+V)
\frac{r}{x(1-x)}\frac{dx}{dr}.
\end{align*}
In fact estimates of the commutators
\begin{align*}
\|[\triangle, \mathcal{A}]\phi\|_n&\leq C(|\vec{a}|_2\|\phi\|_{n+3}+|\vec{a}|_{n+5}\|\phi\|_0),         \\
\|[\triangle, \mathfrak{a}],\phi\|_n&\leq  C (|\vec{a}|_3\|\phi\|_{n+2}+
|\vec{a}|_{n+5}\|\phi\|_0),      \\
\|[\triangle,J]\phi\|_n&\leq C(|\vec{a}|_4\|\phi\|_{n+1}+
|\vec{a}|_{n+5}\|\phi\|_0 )  
\end{align*}
can be derived as in \cite{MakinoFE} and used to prove the elliptic
a priori estimate by induction on $n$.

On the other hand the energy estimate should read
$$ \|\vec{H}\|\leq
C\Big(\|\vec{H}|_{t=0}\|+\int_0^T
\|\vec{G}(t')\|dt'\Big),
$$
where
$$\|\vec{H}\|=(\|H\|^2+\|\dot{D}H\|^2+\|K\|^2)^{1/2} \quad \mbox{with}
\quad \|\cdot\|=\|\cdot\|_{L^2(x^{\frac{N}{2}-1}dx)}, $$
for any solution $\vec{H}=(H,K)^T$ of 
$$\frac{\partial\vec{H}}{\partial t}+\mathfrak{A}\vec{H}=\vec{G},
\quad H|_{x=1}=0, $$
which may not vanish at $t=0$,
so that
$$\|\partial_t^n\vec{h}\|\leq C\Big(
\|\partial_t^n\vec{h}|_{t=0}\|+
\int_0^t
\|\partial_t^n\vec{g}\|+
\int_0^t\|[\partial_t^n,\mathfrak{A}]\vec{h}\|\Big).
$$

Moreover we have an estimate
$$\|\partial_t^{n+1}\vec{h}|_{t=0}\|\leq C(1+
W_n(\vec{g})+
|\vec{a}|_{n+3}^{\langle 0\rangle}),
$$
where
$$W_n(\vec{g})=\sum_{j+k\leq n}\|\partial_t^j\vec{g}|_{t=0}\|_k,$$
provided that $|\vec{a}|_4$ and $W_0(\vec{g})$ are bounded. In order to verify this, it is sufficient to show
$$\|\partial_t^{n+1}\vec{h}|_{t=0}\|_k\leq
C(W_{n+k}(\vec{g})+
|\vec{a}|_{n+k+3}^{\langle 0\rangle}
W_0(\vec{g})+
|\vec{a}|_{k+4}^{\langle 0\rangle}
W_{n-1}(\vec{g}))$$
inductively on $n$ using
$$\|\mathfrak{A}\vec{h}\|_n\leq C
(\|\vec{h}\|_{n+1}+|\vec{a}|_{n+4}
\|\vec{h}\|_0).
$$\\

Then the same discussion using the auxiliary quantity
$$Z_n(\vec{h})=\sum_{j+k=n}\|\partial_t^j\vec{h}\|_k$$
as \cite{MakinoFE} leads us to the estimate
$$
\|\vec{h}\|_{n+1}^{\langle t\rangle}\leq
C\Big(1+
\int_0^t\|\vec{g}\|_{n+1}^{\langle t'\rangle}dt'+
W_n(\vec{g})
+\|\vec{g}\|_n^{\langle T\rangle}+|\vec{a}|_{n+4}^{\langle T\rangle}
\Big)
$$
for $0\leq t\leq T$. 

This estimate for the split problem is sufficient to get the tame estimate for the original $\vec{h}=\vec{h}^{[0]}+\vec{h}^{[1]}$ as \cite{MakinoOJM}. We omit the repetition of the discussion.

\end{document}